\documentclass[a4paper,10pt,reqno]{amsart}
\usepackage[utf8]{inputenc}
\usepackage{amsmath,amssymb,bm, tikz, cancel, amsthm}

\usepackage{upgreek}
\usepackage[mathscr]{euscript}
\usepackage[colorlinks=true,linkcolor=blue,citecolor=orange]{hyperref}

\usepackage{parskip, color}
\usepackage[mathscr]{euscript}
\usepackage{setspace}

\newtheorem{theorem}{Theorem}[section]
\newtheorem{prop}[theorem]{Proposition}
\newtheorem{lemma}[theorem]{Lemma}

\newtheorem{remark}[theorem]{Remark}

\numberwithin{equation}{section}

\newcommand{\su}{u}
\newcommand{\sv}{\mathtt{v}}
\newcommand{\suu}{\mathtt{u}}

\renewcommand{\paragraph}[1]{{\bf #1}}

\DeclareMathOperator{\eps}{\varepsilon}
\DeclareMathOperator{\dxk}{\mathit k}
\DeclareMathOperator{\Z}{\mathbb{Z}}
\newcommand{\mbf}[1]{\boldsymbol{#1}}
\DeclareMathOperator{\fq}{\omega}
\DeclareMathOperator{\Mpert}{\mathfrak{m}}
\newcommand{\zxk}{z}
\DeclareMathOperator{\Ld}{\mathit{L}}
\DeclareMathOperator{\Nd}{\Ld}
\DeclareMathOperator{\dxp}{\mathit{x}}
\DeclareMathOperator{\idf}{\omega}
\DeclareMathOperator{\idk}{k}
\DeclareMathOperator{\sC}{\mathtt{C}}
\DeclareMathOperator{\Ks}{\varpi_s}

\DeclareMathOperator{\Mpt}{\Delta}

\begin{document}

\author{Josselin Garnier}
\address{Centre de Math\'{e}matiques Appliqu\'{e}es, Ecole Polytechnique, Institut Polytechnique de Paris, 91120 Palaiseau, France}
\email{josselin.garnier@polytechnique.edu}

\author{Basant Lal Sharma}
\address{Department of Mechanical Engineering, Indian Institute of Technology Kanpur, Kanpur, 208016 UP, India}
\email[Corresponding author]{bls@iitk.ac.in}

\title[Effective dynamics in randomly perturbed lattices]{
Effective dynamics in lattices with random mass perturbations}

\maketitle

\date{\today}

\begin{abstract}
We consider a one-dimensional mono-atomic lattice with random perturbations of masses spread over a finite number of particles. Assuming Newtonian dynamics and linear nearest-neighbour interactions and allowing for a provision of pinning due to substrate interaction, we discuss 
a transient dynamics problem and a time-harmonic transmission problem. By a stochastic, multiscale analysis we provide asymptotic expressions for the displacement field that propagates through the random perturbations and for the time-harmonic transmission coefficients. 
These theoretical predictions are supported by illustrations of their agreements with numerical simulations. 
\end{abstract}

\maketitle

\section{Introduction}
In this paper we consider a one-dimensional chain of particles. 
Each particle interacts through a nearest-neighbor potential. 
The difference equation that governs the dynamics of the one-dimensional lattice is deduced
from Newton’s law.

In the time-dependent framework, the problem for the displacement field has the form
\begin{equation}\begin{split}
(1+{{\Mpert}}_{\dxp})\ddot{\su}_{\dxp}(t)
&={\su}_{{\dxp}+1}(t)+{\su}_{{\dxp}-1}(t)-(2+\Ks) {\su}_{\dxp}(t),\quad \dxp\in\Z, t\in\mathbb{R}  ,
\label{TDeq}
\end{split}\end{equation}
where
$\Ks\geq 0$, the masses of the particles are $1+{\Mpert}_{\dxp}$,
and the initial condition is
\begin{equation}\begin{split}\begin{split}
{\su}_{\dxp}(0)={\suu}_{\dxp}^{(0)},\quad \dot{\su}_{\dxp}(0)={\sv}_{\dxp}^{(0)}, \quad\dxp\in\Z,
\label{TDiniteq2}
\end{split}\end{split}\end{equation}
with a specified ${\suu}^{(0)}$ and ${\sv}^{(0)}$ in $\ell_2(\Z)$.
We are particularly interested in solving for $\su:\Z\times \mathbb{R}\to\mathbb{R}$ with the initial condition
\begin{equation}\begin{split}\begin{split}
{\su}_{\dxp_0}(0)=1,{\su}_{\dxp}(0)=0,~~\dxp\in\Z\setminus\{\dxp_0\}, \dxp_0\in\Z,\quad\qquad
\dot{\su}_{\dxp}(0)=0, ~~\dxp\in\Z.
\label{TDiniteq1}
\end{split}\end{split}\end{equation}

\begin{figure}[htb!]
\centering
\includegraphics[width=\textwidth]{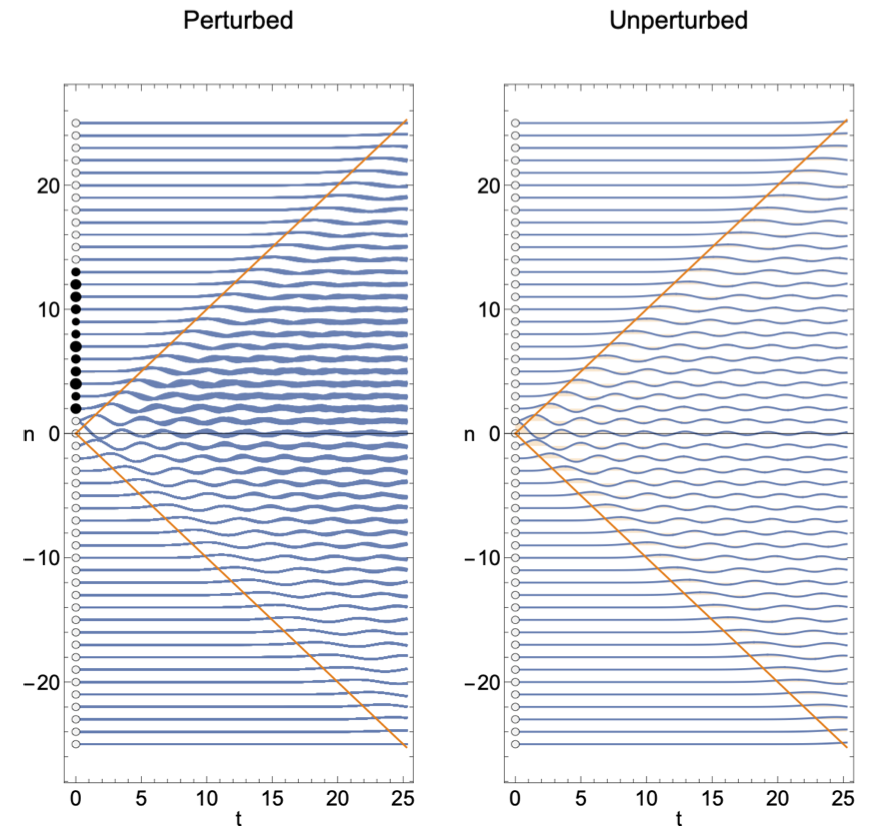}
\caption{Time domain problem with (left) and without (right) mass perturbations. Here ${\Ld}=12, \sigma=0.15$ and $\Ks=0$ (the number of curves on the left is $51$ corresponding to a set of realizations of random mass perturbations with the same statistics). Perturbed masses are indicated as solid black disks on the schematic, adjacent to $n$-axis, whereas empty dots refer to the regular lattice. Mass perturbation are located from $n=2$ to $n=13$. Initial condition \eqref{TDiniteq1} is supported outside the mass defect with $\dxp_0$ corresponding to $n=0$.
The orange lines indicate boundary of cone with unit speed of propagation.}
\label{timedomainonedimlatticerandommasspertb}
\end{figure}

\begin{figure}[htb!]
\centering
\includegraphics[width=\textwidth]{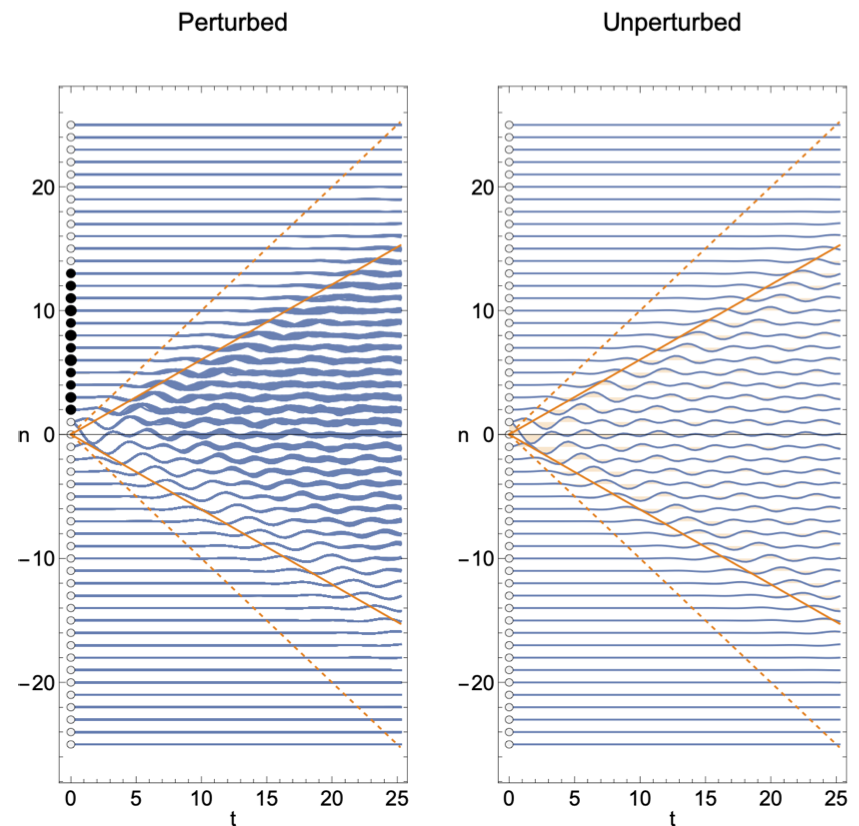}
\caption{Time domain problem with (left) and without (right) mass perturbations. Here $\Ks=1.1$ and the other parameters and details are the same ones as in Fig. \ref{timedomainonedimlatticerandommasspertb}. 
%Perturbed masses are indicated as solid black disks whereas empty dots refer to the unperturbed lattice. Mass perturbations are located from $n=2$ to $n=13$. Initial condition \eqref{TDiniteq1} is supported outside the mass defect with $\dxp_0$ corresponding to $n=0$.
{The dashed orange lines indicate boundary of cone with unit speed of propagation,
the solid orange lines indicate boundary of cone with speed of propagation $1/\alpha_s$ (see Eq.~(\ref{eq:defalphaomegas}))}.}
\label{timedomainonedimlatticerandommasspertbKs}
\end{figure}

Figs. \ref{timedomainonedimlatticerandommasspertb} and \ref{timedomainonedimlatticerandommasspertbKs} present trajectories of particles in the lattice obtained by solving \eqref{TDeq}, \eqref{TDiniteq1} using standard numerical method assuming $\Ks=0$ and $\Ks\ne0$, respectively, ${\Mpert}_{\dxp}$ are independent and identically distributed with mean zero and variance $\sigma^2$ in the section $\dxp \in [1,{\Ld}]\cap\Z$, and ${\Mpert}_{\dxp}=0$ outside the section $[1,{\Ld}]\cap\Z$. One can observe that the mass perturbations induce perturbations in the dynamics that we %will 
describe in Section \ref{sec:td}.

Equation \eqref{TDeq} was studied in \cite{Schrodinger}, where several  remarkable analyses and features (such as a closed-form expression of the solution when $\Ks=0$) were proposed; see also \cite{Muhlich,Charlotte}.
Equation \eqref{TDeq}  describes the vibration of an infinite mono-atomic chain with nearest-neighbour interactions and belongs to a class of problems that appear in the study of the dynamics of crystal lattices \cite{Brillouin}.

\begin{figure}[ht!]
(a){\includegraphics[width=.45\textwidth]{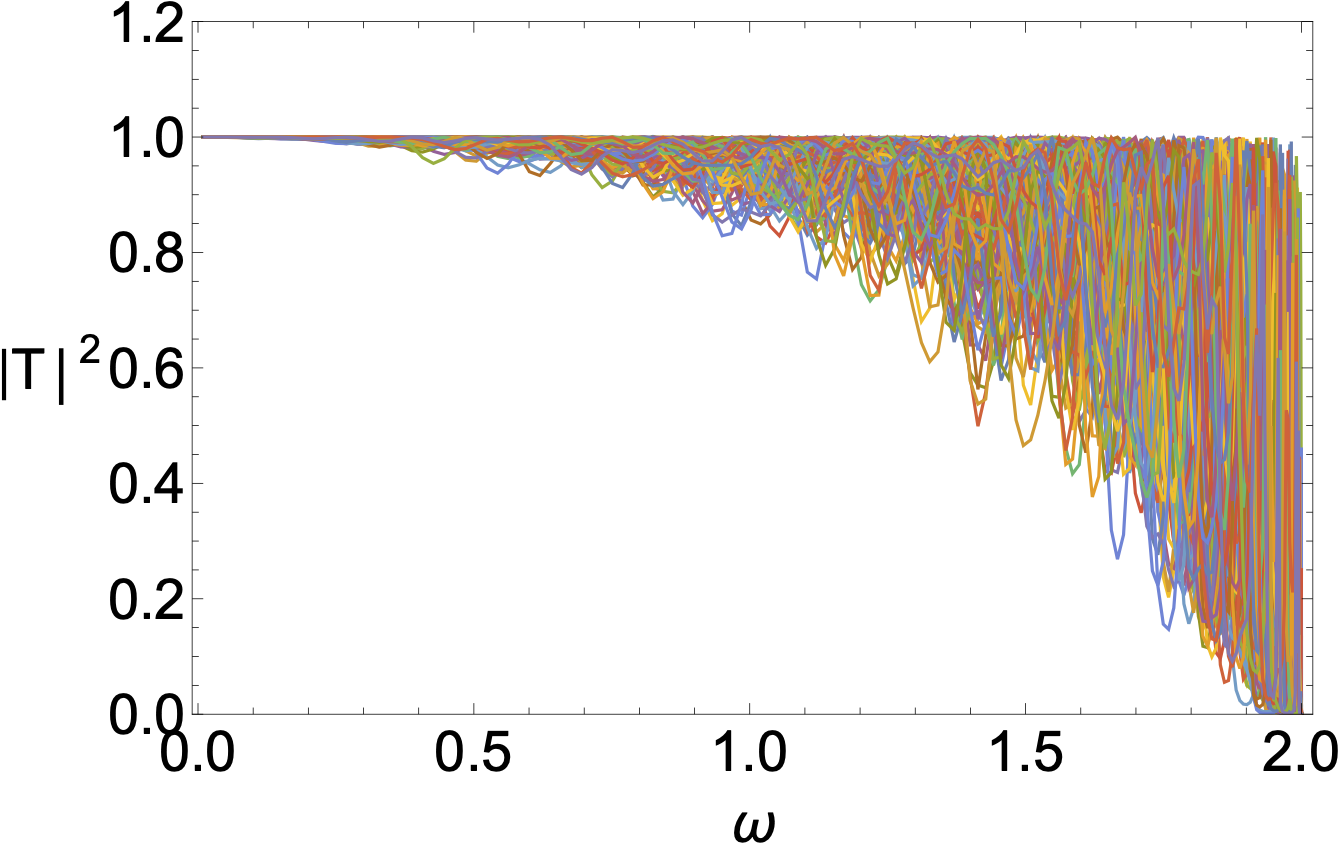}}(b){\includegraphics[width=.45\textwidth]{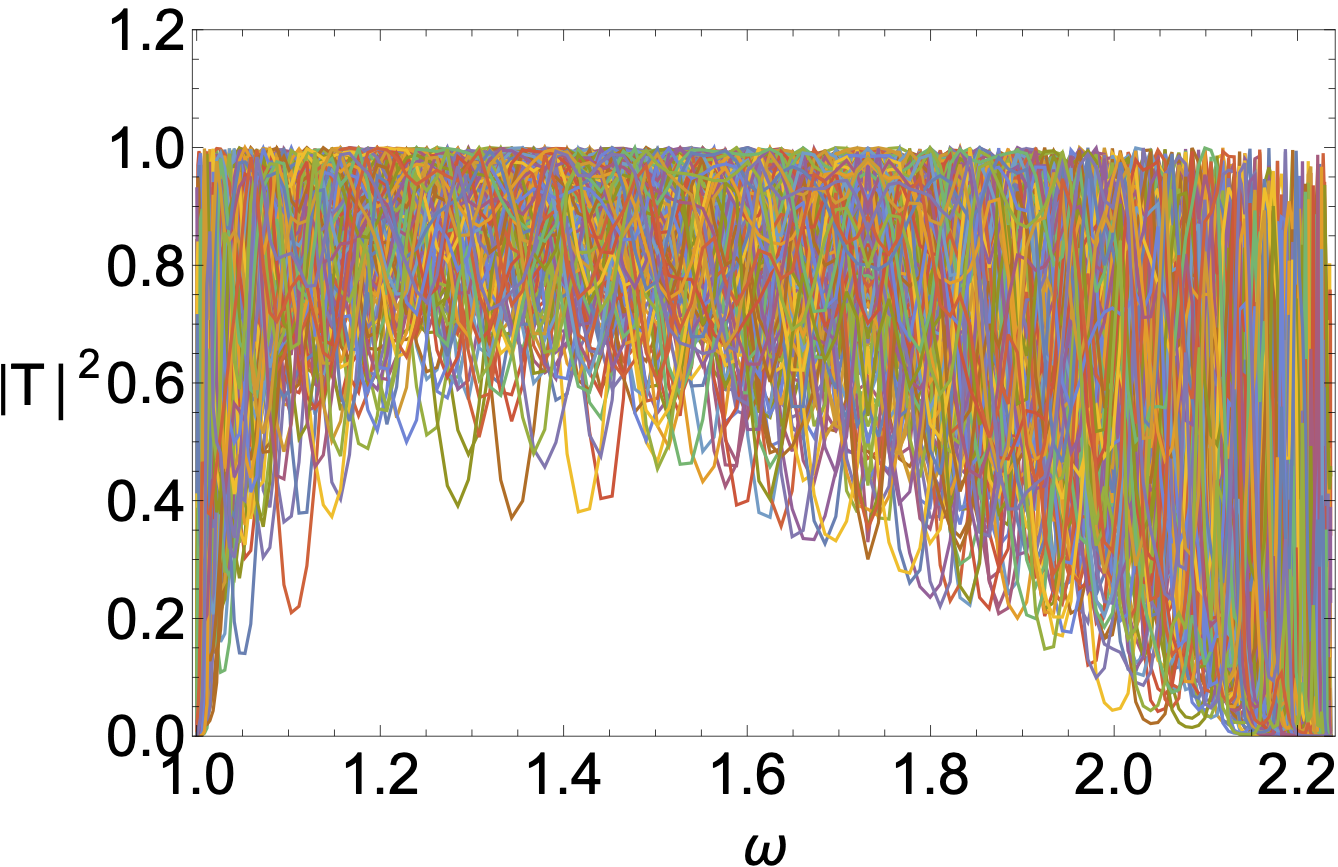}}
\caption{Numerical evaluations of $|T|^2$ using the exact expression \eqref{Tfullmatched} described in Appendix \ref{Tfullmatchedapp} with (a) $\Ks={0}$ (b) $\Ks={1}$, $\sigma=0.05, {\Ld}=40$ and ${200}$ different realizations of the mass perturbations ${{\Mpert}}_{\dxp}$, $ {\dxp} \in [1,{\Ld}]$.
In the absence of mass perturbation on $[1,\Ld]\cap\Z$, $T\equiv1$.}
\label{freqdomainonedimlatticerandommasspertb}
\end{figure}

In the time-harmonic framework with frequency $\omega$ the solution of \eqref{TDeq} has the form
\begin{equation}\begin{split}
\label{eq:thwave}
\su_{\dxp}(t) = e^{- i {\idf} t} \hat{\su}_{\dxp},\quad \dxp\in\Z, t\in\mathbb{R},
\end{split}\end{equation}
where $\hat{\su}: \Z\to\mathbb{C}$ is a solution of 
\begin{equation}\begin{split}
\label{eq:th}
-{\idf}^2 (1+{\Mpert}_{\dxp}) \hat{\su}_{\dxp} = (\hat{\su}_{{\dxp}+1}+\hat{\su}_{{\dxp}-1}-2\hat{\su}_{\dxp}) - \Ks \hat{\su}_{\dxp},\quad \dxp\in\Z.
\end{split}\end{equation}
Fig. \ref{freqdomainonedimlatticerandommasspertb} presents the transmittance obtained by solving \eqref{eq:th} using numerical methods for different realizations of the random mass perturbations ${\Mpert}_{\dxp}$ in the section  $\dxp \in [1,{\Ld}]\cap\Z$. One can observe that the  transmittance has strong fluctuations and we %will 
describe its statistics in Section \ref{sec:th}.

Besides its application to lattice vibrations, the equation \eqref{eq:th} also belongs to a class of discrete scattering problems in the context of the discrete Schr\"{o}dinger equation (within tight-binding model of the electrons in crystals) \cite{Bloch,Slater}; it has also played a crucial role in the discovery of significant phenomena such as the famous localization result of \cite{Anderson}. 
In the domain of electrical engineering, network synthesis and filter design \cite{Kuo}, LC circuits based lattice structures also involve similar difference operators as in \eqref{TDeq} and \eqref{eq:th}, while such operators also appear in the lumped circuit models for electromagnetic metamaterials \cite{Itoh,Deymier}.
The equation \eqref{eq:th} naturally emerges in case of time harmonic lattice waves in one dimension \cite{Brillouin,BornHuang1912,BornHuang1985}.

% need to compress the sentences below in next three paragraphs
In fact, over the last century till this date, the analysis of one-dimensional lattice models accounting for disorder and randomness, has been a part of several physics-oriented and mathematical researches \cite{Han}.

In the backdrop of quantum mechanics of electronic wave functions, the theorem due to \cite{Bloch} was utilized by \cite{Kronig} leading to a `binding' introduced by the potential field in a one-dimensional problem. The same model of a one-dimensional crystal was studied further by \cite{Saxon} who used the scattering-matrix method to study the energy levels and wave functions for the problems of a single impurity. In particular, for linear chains, but now in framework of phonons, the normal-mode frequency spectra of certain binary isotropically disordered harmonic lattices were presented by \cite{Payton}. The study of disordered systems was reviewed by \cite{Elliott} during these developments;
several elementary excitations were seen through one common descriptive Hamiltonian and similarity between these problems was described including the case of random material heterogeneity. 
A critical survey of the literature on electronic transmission models of one-dimensional disordered solids was presented by \cite{Erdos} who also connected it to the problem of elastic vibrations in disordered solids. With advent of computers,
some calculations for narrow wires with disorder following the Anderson model \cite{Anderson} were presented by \cite{Godin88}. The study on this class of problems continued for several decades during the increasing dominance of numerical simulations in research. 

Over last few years, there have been some exciting developments in a similar set of research problems too.
According to \cite{Chaudhuri}, the disorder causes short-wavelength phonon modes to be localized so the heat current in this system is carried by the extended phonon modes which can be either diffusive or ballistic; \cite{Dhar} has reviewed heat conduction in harmonic crystals in connection with the Landauer formula for phonon heat transfer.
\cite{Komorowski,Komorowski2} considered the long time limit for the solutions of a discrete wave equation with weak stochastic forcing. 
\cite{Basile2014} proved a diffusive behaviour of the energy fluctuations in a system of harmonic oscillators with a stochastic perturbation of the dynamics that conserves energy and momentum; see also \cite{Basile}.
\cite{Jara} considered the unpinned case, as well as the pinned case, finding relations to fractional diffusion equation or classical diffusion while \cite{Tomasz} 
incorporated Ornstein-Uhlenbeck process.
\cite{Basile2} also reviewed these rigorous results about the macroscopic behaviour of harmonic chains with the dynamics perturbed by a random exchange of velocities between nearest neighbor particles and the role of fractional heat equation concerned with pinned or unpinned systems in any dimension. 
\cite{Sivak} studied the heat fluctuations in a one-dimensional finite harmonic chain of interacting active Ornstein-Uhlenbeck particles with the chain ends connected to heat baths.
Indeed, impurities and vacancies play an important role in the thermal conductivity of materials specially due to local mass distortions around point defects; See \cite{Amir,Xie,Lindsay, Katcho, Dongre,Thebaud} as examples of some recent studies from physics related viewpoint involving mass disorder in crystalline materials.

In the simple physical model that we study in this article, the mass perturbation is random, but time-independent and spread over a finite number of particles; the two semi-infinite ends of lattice carry the phonon modes transmitted across such perturbation. In a sense, we initiate the development of certain stochastic framework for capturing effective dynamical behaviours of discrete random media by using a prototype of one dimensional lattice model in this article. 
The general questions posed in the present article concerning \eqref{TDeq} and \eqref{eq:th} are also related to quasi-one-dimensional problems of waveguides; see \cite{Ong}, for example. For example, the problem of electronic transport through several interfaces and junctions, naturally described in terms of transmission and reflection of each interface, following the Landauer-B\"{u}ttiker approach \cite{Landauer1,Buttiker1}. The solution of scattering problem decides the conductance of such single interface in the linear response regime; see, for example, \cite{Blsbif,Blshexa,Blshexa2,Blsstep} for an application to transport in waveguides where the forward problem of quasi-one-dimensional discrete Schr\"{o}dinger equation has been solved exactly in case of single interface but a non-compact perturbation. On the other hand, the question of inverse problem is also related to such forward analysis of discrete media and attaining information about statistics of the perturbation can be useful. A recent result on inverse scattering on lattices, in particular in one dimension, is given in \cite{novikovbls}; also see the references therein for the general problem of inverse scattering in discrete framework. 
% need to compress the above sentences

The paper is organized as follows. In Section \ref{sec:prel} the unperturbed problem is addressed and asymptotic solutions are presented. The time-dependent, resp. time-harmonic, problem with random mass perturbations is studied in Section \ref{sec:td}, resp. \ref{sec:th}.
Proofs are given in the following sections.

\section{Preliminary results}
\label{sec:prel}
In this section we consider the unperturbed case ${{\Mpert}}_{\dxp}=0$ for all $x \in \mathbb{Z}$.
The following lemmas are proved in Section~\ref{sec:proofprel}.

\begin{lemma}[Solution without perturbation for {$\Ks=0$}{}]
The solution of (\ref{TDeq}) with the initial conditions (\ref{TDiniteq1}) in the absence of perturbations ${\Mpert}_{\dxp} \equiv 0$ and with ${\Ks}=0$ is of the form
\begin{equation}
\label{eq:relieuhatuhomo1}
{\su}_{{\dxp}}(t) = \frac{1}{2\pi} \int_0^{+\infty} \hat{\su}_{{\dxp}}(\fq) e^{-i \fq t} d\fq + c.c., \qquad t \geq 0, \quad \dxp \in \Z,
\end{equation}
where the Fourier components are given by 
\begin{equation}
\label{eq:relieuhatuhomo2}
\hat{\su}_{{\dxp}}(\fq) =  
\hat{c}(\fq) \cos \big( k(\fq) (\dxp-\dxp_0) \big) , 
\end{equation}
with
\begin{equation}
\label{eq:defkc}
k(\fq) = 
2 \, {\rm arcsin}\big(\frac{\fq}{2}\big),\qquad 
\hat{c}(\fq) = \frac{2}{\sqrt{4-\fq^2}} {\bf 1}_{(0,2)}(|\fq|)  .
\end{equation}
\label{lemma1}
\end{lemma}

We can also write since $\hat{\su}_{{\dxp}}(-\fq) =\overline{ \hat{\su}_{{\dxp}}(\fq) }$ that
\begin{equation}
{\su}_{{\dxp}}(t) = \frac{1}{2\pi} \int_{-\infty}^{+\infty} \hat{\su}_{{\dxp}}(\fq) e^{-i \fq t} d\fq  , \qquad t \geq 0, \quad \dxp \in \Z. 
\label{eq:expressuy2homo}
\end{equation}

\begin{remark}
We can make the change of variable $\fq \mapsto 2 \sin s$ in (\ref{eq:expressuy2homo}) and we obtain
\begin{equation}\begin{split}
{\su}_{{\dxp}}(t) = \frac{1}{2\pi} \int_{-\pi/2}^{\pi/2} \big( e^{2i s (\dxp-\dxp_0)} + e^{-2i s(\dxp-\dxp_0) }\big) e^{-2i  \sin(s) t} ds , \qquad t \geq 0, \quad \dxp \in \Z,
\end{split}\end{equation}
which gives 
\begin{equation}\begin{split}
{\su}_{{\dxp}}(t) = J_{2(\dxp-\dxp_0)}(2t),\qquad t \geq 0, \quad \dxp \in \Z,
\end{split}\end{equation}
where $J$ is the Bessel function of the first kind. This expression is valid only for $\Ks=0$ and there is no equivalent expression for $\Ks>0$.
\end{remark}

We can study the asymptotic behavior of $ {\su}_{{\dxp_0+\dxp}}(t)$ for large $t$ and $\dxp$.
Let us consider $\dxp \gg1 $ and $t = \alpha \dxp$, $\alpha \in (1,+\infty)$. 
A stationary phase argument gives the following result.
\begin{lemma}[asymptotic behavior $\Ks=0$]
For any $\alpha \in (1,+\infty), \dxp \in \Z$, we have 
\begin{equation}
 {\su}_{{\dxp_0+\dxp}}( \alpha \dxp )  = \frac{1}{\sqrt{\pi \dxp} \sqrt[4]{\alpha^2-1}}
 \cos \Big( \frac{\pi}{4}+ 2\big({\rm arccos}(1/\alpha) - \sqrt{\alpha^2-1}\big) \dxp \Big) +o(\frac{1}{\sqrt{\dxp}}),
 \label{eq:asyhomouy1}
\end{equation}
as $\dxp \to +\infty$.
\label{lemma2}
\end{lemma}
For $\alpha \in [0,1)$ there is no stationary point which implies that $ {\su}_{{\dxp_0+\dxp}}( \alpha \dxp ) $ is much smaller than ${1}/{\sqrt{\dxp}}$.
For $\alpha=1$, i.e. at times close to $\dxp$, a refined study shows the behavior of the front field.
\begin{lemma}[asymptotic front behavior $\Ks=0$]
For any $\beta \in \mathbb{R}, \dxp \in \Z,$ we have 
\begin{equation}
 {\su}_{{\dxp_0+\dxp}}(\dxp+\sqrt[3]{\dxp} \beta)  = \frac{1}{ \sqrt[3]{\dxp}} 
 Ai(-2\beta) +o(\frac{1}{\sqrt[3]{\dxp}})   ,
 \label{eq:asyhomouy2}
 \end{equation}
as $\dxp \to +\infty$,
where $Ai$ is the Airy function.
\label{lemma3}
\end{lemma}

\begin{remark}
Recall that $Ai({\dxp}) \sim e^{ -({2}/{3}) {\dxp}^{3/2}} /(2 \sqrt{\pi}\sqrt[4]{{\dxp}})$ for ${\dxp} \gg 1$ which shows that the field ${\su}_{{\dxp_0+\dxp}}(t)$ is vanishing before the arrival time $\dxp$ and the field is essentially contained in the cone with unit speed of propagation  as seen in Figure~\ref{timedomainonedimlatticerandommasspertb} right.
Moreover $Ai(-{\dxp}) \sim\cos [  ({\pi}/{4})- ({2}/{3}) {\dxp}^{3/2}] /(\sqrt{\pi} \sqrt[4]{{\dxp}})$ for $ {\dxp} \gg 1$, which shows that the two expansions (\ref{eq:asyhomouy1}) and (\ref{eq:asyhomouy2}) match.
\end{remark}

If ${\Ks}>0$, then the same procedure gives the expression (\ref{eq:relieuhatuhomo1}-\ref{eq:relieuhatuhomo2}) for the solution ${\su}_{{\dxp}}(t) $ with modified expressions for $k(\fq)$ and $\hat{c}(\fq)$.

\begin{lemma}[Solution without perturbation for {${\Ks}>0$}{}]
If ${\Ks}>0$, then the solution ${\su}_{{\dxp}}(t) $ has the form (\ref{eq:relieuhatuhomo1}-\ref{eq:relieuhatuhomo2}) with
\begin{equation}
\label{eq:defkomegas}
k(\fq) = 
2 \, {\rm arcsin}\big(\frac{\sqrt{\fq^2-{\Ks}}}{2}\big),\qquad 
\hat{c}(\fq) = \frac{2|\fq|}{
{\sqrt{\fq^2-{\Ks}}\sqrt{4+{\Ks} -\fq^2 }}
} {\bf 1}_{(\sqrt{{\Ks}},\sqrt{{\Ks}+4})}(|\fq|) .
\end{equation}
\end{lemma}

We can study the asymptotic behavior of $ {\su}_{{\dxp_0+\dxp}}(t)$ for large $t>0$ and $\dxp \in \Z$.
Let us consider $\dxp \gg 1 $ and $t = \alpha \dxp$, $\alpha \in (\alpha_s,+\infty)$, with
\begin{equation}
\label{eq:defalphaomegas}
\alpha_s=\frac{\sqrt{2}}{\sqrt{2+{\Ks}-\fq_s^2}},\qquad \fq_s = \sqrt[4]{4{\Ks}+{\Ks}^2} .
\end{equation}
Note that $\alpha_s$ is larger than one.
A stationary phase argument gives the following result.
\begin{lemma}[asymptotic behavior ${\Ks}>0$]
For any $\alpha \in (\alpha_s,+\infty), \dxp \in \Z,$ we have
\begin{align}
\nonumber {\su}_{{\dxp_0+\dxp}}( \alpha \dxp )  = & \frac{ \sqrt{2} {\fq_{\alpha}^{+}}^{3/2} }{\sqrt{\pi \alpha ( {\fq_{\alpha}^{+}}^4-\fq_s^4) \dxp}}
\cos\Big( \frac{\pi}{4} +  [k(\fq_{\alpha}^{+})  -\fq_{\alpha}^{+} \alpha] \dxp\Big)
\\
& + \frac{\sqrt{2} {\fq_{\alpha}^{-}}^{3/2} }{\sqrt{\pi \alpha (\fq_s^4 -  {\fq_{\alpha}^{-}}^4) \dxp}}
\cos\Big( -\frac{\pi}{4} +  [k(\fq_{\alpha}^{-})  -\fq_{\alpha}^{-} \alpha] \dxp\Big)
+o(\frac{1}{\sqrt{\dxp}}),
 \label{eq:asyhomouy1:ks}
\end{align}
as $\dxp \to +\infty$,
where
\begin{equation}
\label{eq:defomagealphaks}
{\fq_{\alpha}^{\pm}}^2 =  \frac{2}{\alpha^2} \Big[ (1+\frac{{\Ks}}{2})\alpha^2 -1
\pm \sqrt{\alpha^4 -(2+{\Ks})\alpha^2 +1}\Big].
\end{equation}
\label{lemma4}
\end{lemma}

If $\alpha \in [0,\alpha_s)$ then there is no stationary point so we can conclude that ${\su}_{{\dxp_0+\dxp}}( \alpha \dxp )$ is much smaller than ${1}/{\sqrt{\dxp}}$.
For $\alpha=\alpha_s$, i.e.
at times close to $\alpha_s \dxp$, a refined study gives the behavior of the front field.
\begin{lemma}[asymptotic front behavior ${\Ks}>0$]
For any $\beta \in \mathbb{R}, \dxp \in \Z,$ we have 
\begin{equation}\begin{split}
 {\su}_{{\dxp_0+\dxp}}(\alpha_s \dxp+\beta \sqrt[3]{\dxp})  &= \frac{\sqrt[3]{2}}{\sqrt[3]{\dxp}} 
 Ai\Big(- \sqrt[3]{2} \frac{\beta}{\alpha_s}\Big)
 \cos\Big( [k(\fq_s) - \fq_s \alpha_s] \dxp -  \fq_s \beta \sqrt[3]{\dxp}\Big)\\
 & +o(\frac{1}{\sqrt[3]{\dxp}})
 , 
 \label{eq:asyhomouy2:ks}
\end{split}\end{equation}
as $\dxp \to +\infty$.
\label{lemma5}
\end{lemma}

Since $Ai({\dxp})$ decays very fast for ${\dxp} \gg 1$, this confirms that the field is vanishing before the arrival time $\alpha_s \dxp$.
The field is essentially contained in the cone with speed of propagation $1/\alpha_s$ as seen in Figure~\ref{timedomainonedimlatticerandommasspertbKs} right.

\section{Effective dynamics for the time-dependent problem}
\label{sec:td}
In this section we assume that, for $\dxp \in [1,L]\cap\Z$, the variables ${\Mpert}_{\dxp}$ are independent and identically distributed with mean zero and variance $\sigma^2$:
\begin{equation}\begin{split}
\mathbb{E}[{\Mpert}_{\dxp}^2]=\sigma^2.
\end{split}\end{equation}
The forthcoming results are obtained by a multiscale analysis that is valid when 
$\sigma \ll 1$ and ${\Ld}$ is of the order of $\sigma^{-2}$ and they are proved in Section  \ref{sec:prooftd}. 
We consider the initial condition \eqref{TDiniteq1} with $\dxp_0=0$.
We define
\begin{equation}\begin{split}
\label{def:gamma:ind}
\gamma(\fq)= \frac{\sigma^2{\fq}^4 }{4 \sin^2 k(\fq) },
\end{split}\end{equation}
where $k$ is given by (\ref{eq:defkomegas}).

\begin{theorem}[Mean field with random perturbations and {$\Ks=0$}]
If $\sigma \ll 1$, $\dxp > {\Ld} \gg 1, \dxp \in \Z$, then
the mean field of \eqref{TDeq}, \eqref{TDiniteq1} has the form for $\alpha \in (1,+\infty)$:
\begin{equation}\begin{split}
\mathbb{E} \big[ {\su}_{{\dxp}}( \alpha \dxp ) \big] = \frac{1}{\sqrt{\pi \dxp} \sqrt[4]{\alpha^2-1}}
\cos \Big( \frac{\pi}{4}+ 2\big({\rm arccos}(1/\alpha) - \sqrt{\alpha^2-1}\big) \dxp \Big) \\
\times e^{- \gamma({\fq}_\alpha) {\Ld}} +o(\frac{1}{\sqrt{\dxp}}),
\label{eq:espasyranduy1}
\end{split}\end{equation}
where 
\begin{align} 
\label{def:omegaalpha}
{\fq}_\alpha &= \frac{2\sqrt{\alpha^2-1}}{\alpha} ,
\\
\gamma({\fq}_\alpha) &= \sigma^2 (\alpha^2-1).
\end{align}
\label{thmTDmeanks0}
\end{theorem}

\begin{theorem}[Mean front with random perturbations and {$\Ks=0$}]
At times close to $\dxp>0, \dxp\in \Z,$ we have
\begin{equation}\begin{split}
\mathbb{E} \big[ {\su}_{{ \dxp}}(\dxp+\sqrt[3]{\dxp} \beta) \big] = \frac{1}{ \sqrt[3]{\dxp}} 
Ai(-2\beta) +o(\frac{1}{\sqrt[3]{\dxp}})
,\quad \beta \in \mathbb{R}.
\label{eq:espasyranduy2}
\end{split}\end{equation}
\label{thmTDfrontks0}
\end{theorem}

The expressions (\ref{eq:espasyranduy1}-\ref{eq:espasyranduy2}) show that the mean field for times $t$ close to the front $\dxp$ is not affected to leading order by the random perturbations, but the mean field following the front for times $t $ larger than $ \dxp $ is affected.
This is in contrast with the results known for the scalar wave equation in which the wave front experiences two different phenomena: a deterministic attenuation and spreading and a random time shift. The attenuation and spreading is described by a deterministic kernel determined by the statistics of the random medium. The random time shift has Gaussian statistics with mean zero and variance that depends on the statistics of the random medium.
The stabilization of the wave front in randomly layered media was first noted by O’Doherty and Anstey in a geophysical context \cite{oda}. A time-domain integral equation approach was given in \cite{white88,burridge89}. A frequency-domain approach was presented in \cite{clouet94,book}.

\begin{remark}[Transmitted field]
In the region $\dxp>{\Ld}, \dxp \in \Z$, 
the solution has the form 
\begin{equation}\begin{split}
{\su}_{{\dxp}}( \alpha \dxp ) = \frac{1}{\sqrt{\pi \dxp} \sqrt[4]{\alpha^2-1}}
\cos \Big( \frac{\pi}{4}+ 2\big({\rm arccos}(1/\alpha) - \sqrt{\alpha^2-1}\big) \dxp +\sqrt{\gamma({\fq}_\alpha)} W_{\Ld} \Big) \\
\times e^{- \frac{\gamma({\fq}_\alpha)}{2} {\Ld}} +o(\frac{1}{\sqrt{\dxp}}),
\label{eq:asyranduy1}
\end{split}\end{equation}
where $W_{\Ld} \sim {\mathcal N}(0,{\Ld})$.
At times close to $\dxp>0, \dxp \in \Z,$ we have
\begin{equation}\begin{split}
{\su}_{{ \dxp}}(\dxp+\sqrt[3]{\dxp} \beta) = \frac{1}{ \sqrt[3]{\dxp}} 
Ai(-2\beta) +o(\frac{1}{\sqrt[3]{\dxp}})
,\quad \beta \in \mathbb{R}.
\label{eq:asyranduy2}
\end{split}\end{equation}
These expressions show that the field for times $t$ close to the front $\dxp$ is not affected to leading order by the random perturbations, but the field following the front for times $t $ larger than $ \dxp $ is affected.
\end{remark}

We now address the case when $\Ks>0$.
Let $\alpha_s$ 
be defined by \eqref{eq:defalphaomegas}.

\begin{theorem}[Mean  field with random perturbations and {$\Ks>0$}]
For $\sigma \ll 1$, $\dxp > {\Ld}\gg 1, \dxp \in \Z$, 
the mean field has the form for $\alpha \in (\alpha_s,+\infty)$:
\begin{align}
\nonumber \mathbb{E}\big[ {\su}_{{\dxp_0+\dxp}}( \alpha \dxp ) \big] = & \frac{ \sqrt{2} {{\fq}^{+}_{\alpha}}^{3/2} }{\sqrt{\pi \alpha ({{\fq}^{+}_{\alpha}}^4-{\fq}_s^4) \dxp}}
\cos\Big( \frac{\pi}{4} + [k({\fq}^{+}_{\alpha}) -{\fq}^{+}_{\alpha} \alpha] \dxp\Big) \notag\\
&\times e^{-\gamma({\fq}^{+}_{\alpha}){\Ld}}\notag
\\
& + \frac{\sqrt{2} {{\fq}^{-}_{\alpha}}^{3/2} }{\sqrt{\pi \alpha ({\fq}_s^4 - {{\fq}^{-}_{\alpha}}^4) \dxp}}
\cos\Big( -\frac{\pi}{4} + [k({\fq}^{-}_{\alpha}) -{\fq}^{-}_{\alpha} \alpha] \dxp\Big)\notag\\
&\times e^{- \gamma({\fq}^{-}_{\alpha}){\Ld}}
+o(\frac{1}{\sqrt{\dxp}}),
\label{eq:asyranduy1:ks}
\end{align}
where $k(\fq)$ is given by (\ref{eq:defkomegas}), ${\fq}_s$ is defined by
(\ref{eq:defalphaomegas}), ${{\fq}^{\pm}_{\alpha}}$ is defined by
(\ref{eq:defomagealphaks}) and
\begin{equation}\begin{split}
\gamma({\fq}^{\pm}_{\alpha}) = \frac{\sigma^2 \alpha^2 {{\fq}^{\pm}_{\alpha}}^2}{4}.
\end{split}\end{equation}
\label{thmTDmeanks}
\end{theorem}

\begin{theorem}[Mean front with random perturbations and {$\Ks>0$}]
At times close to $\alpha_s \dxp, \dxp>0, \dxp \in \Z,$ we have
\begin{equation}\begin{split}
\mathbb{E} \big[ {\su}_{{ \dxp}}(\alpha_s \dxp+\sqrt[3]{\dxp} \beta) \big] = \frac{\sqrt[3]{2}}{\sqrt[3]{\dxp}} 
Ai\Big(-\sqrt[3]{2} \frac{\beta}{\alpha_s}\Big)
\cos\Big( [k({\fq}_s) - {\fq}_s \alpha_s] \dxp - {\fq}_s \beta \sqrt[3]{\dxp}\Big)\\
\times e^{-\gamma ({\fq}_s) {\Ld}}+o(\frac{1}{\sqrt[3]{\dxp}}) 
,\quad \beta \in \mathbb{R},
\label{eq:espasyranduy2:ks}
\end{split}\end{equation}
where $k(\fq)$ is defined by (\ref{eq:defkomegas}), ${\fq}_s$ is defined by
(\ref{eq:defalphaomegas}), and
\begin{equation}\begin{split}
\gamma({\fq}_s) = \frac{\sigma^2 \alpha_s^2 {\fq}_s^2}{4} = \frac{\sigma^2 \sqrt{\Ks}\sqrt{4+\Ks}}{(\sqrt{4+\Ks}-\sqrt{\Ks})^2}.
\label{eq:thmTDfrontks}
\end{split}\end{equation}
\label{thmTDfrontks}
\end{theorem}

As $\Ks\to0^+$, it is noted that the expression \eqref{eq:espasyranduy2:ks} reduces to earlier result for $\Ks=0$ and the attenuation drops to zero. Thus the case $\Ks>0$ is characterized by an attenuation of the front field  for times close to $\alpha_s \dxp$, which is different from the behaviour in the case $\Ks=0$.

\begin{figure}[htb!]
(a)\includegraphics[width=.79\textwidth]{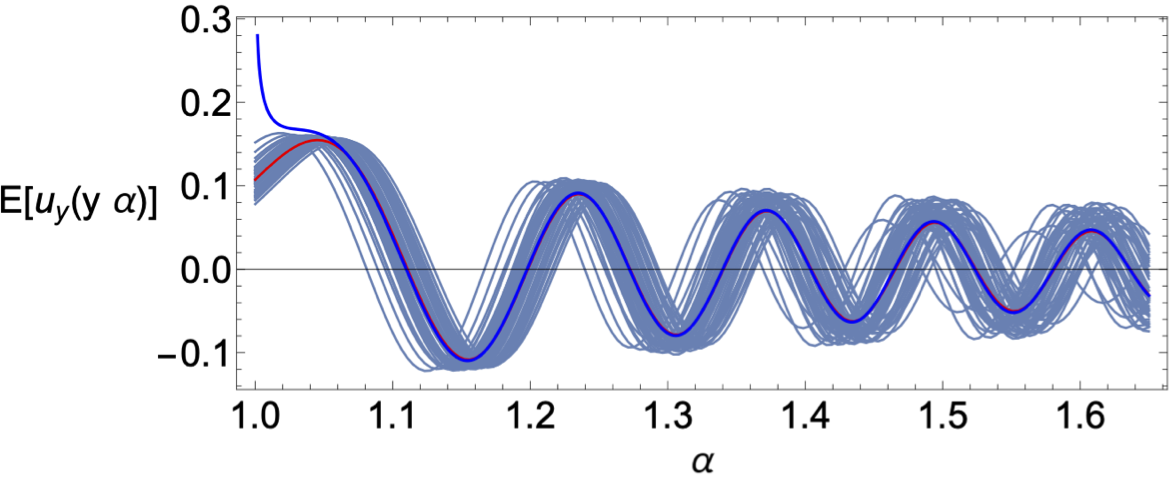}\\
(b)\includegraphics[width=.8\textwidth]{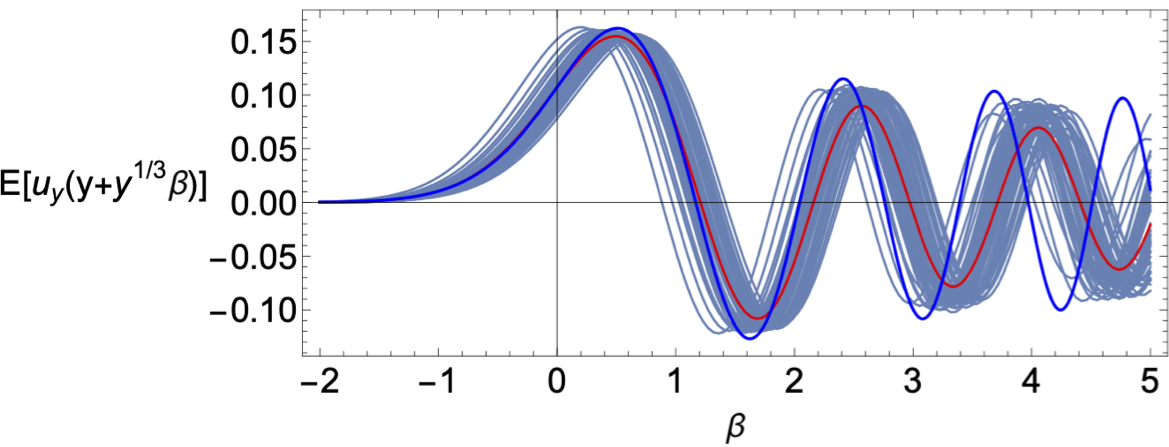}\\
(c)\includegraphics[width=.75\textwidth]{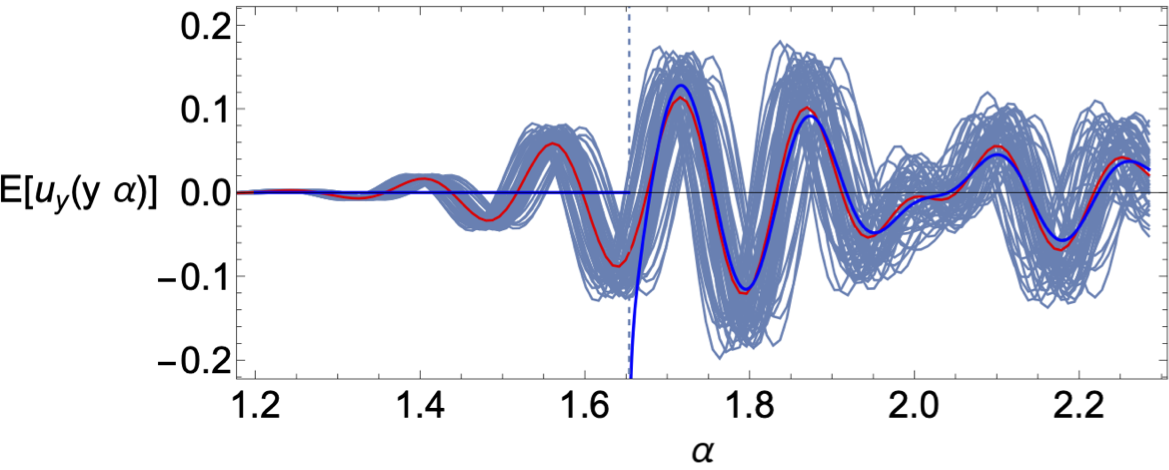}\\
(d)\includegraphics[width=.8\textwidth]{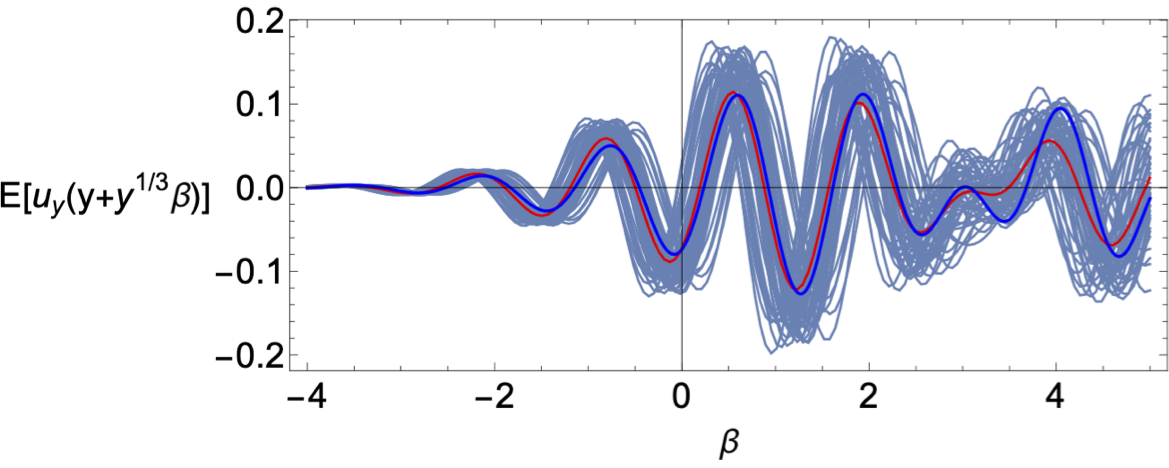}\\
\caption{Time domain problem with random mass perturbations with $\sigma=0.15$, comparing ensemble of numerical solutions (in grayish blue shade with empirical mean curve in red) with asymptotic formulas  (in blue). (a) Mean field \eqref{eq:espasyranduy1}, $\Ks=0$, ${\Nd}=16$. (b) Mean front \eqref{eq:espasyranduy2}, $\Ks=0$, ${\Nd}=16$. (c) Mean field \eqref{eq:asyranduy1:ks}, $\Ks=1.1$, ${\Nd}=8$. (d) Mean front \eqref{eq:espasyranduy2:ks}, $\Ks=1.1$, ${\Nd}=8$.}
\label{TDasymptotics12}
\end{figure}

In Figure \ref{TDasymptotics12} we compare the empirical averages of numerical simulations with the theoretical predictions  for the mean field and the mean front (i.e. \eqref{eq:espasyranduy1}, \eqref{eq:espasyranduy2} when $\Ks=0$ and  \eqref{eq:asyranduy1:ks}, \eqref{eq:thmTDfrontks} when $\Ks>0$).
We obtain  excellent agreement which demonstrates the accuracy of the asymptotic approach.

\section{Effective dynamics for the time-harmonic problem}
\label{sec:th}
Here we assume that the variables ${\Mpert}_{\dxp}$ in the section $\dxp \in [1,{\Ld}]\cap\Z$ are identically distributed with mean zero and variance $\sigma^2$:
\begin{equation}\begin{split}
\mathbb{E}[{\Mpert}_{\dxp}^2]=\sigma^2,
\end{split}\end{equation}
with $\sigma \ll 1$, and ${\Ld}$ is of the order of $\sigma^{-2}$. 

We consider a time-harmonic wave (\ref{eq:thwave}) with $\idf \in (\sqrt{\Ks},\sqrt{\Ks+4})$ (propagative regime).
When $\Mpert_{\dxp} =0$ for ${\dxp} \leq 0, \dxp \in \Z,$ and for ${\dxp}>{\Ld}, \dxp \in \Z,$
and a unit-amplitude right-going input wave is incoming from the left, the solution has the form
\begin{align}
\hat{u}_{\dxp} &= e^{ik{\dxp}}+R e^{-ik{\dxp}} \quad \mbox{ for } {\dxp} \leq 0, \dxp \in \Z,\\
\hat{u}_{\dxp} &= T e^{ik{\dxp}}\quad  \mbox{ for } {\dxp} >{\Ld}, \dxp \in \Z,
\end{align}
where $k(\idf)$ is the solution of the dispersion relation
\begin{equation}
    -\idf^2 =2\cos k-2 -\Ks ,
\end{equation} 
which has the form (\ref{eq:defkomegas}).
The time-harmonic field $\hat{u}_{\dxp}$ satisfies (\ref{eq:th}) for ${\dxp}\in[1,{\Ld}]\cap\Z$.
The complex coefficient $R$, resp. $T$, is the reflection, resp. transmission, coefficient of the perturbed section $[1,{\Ld}]\cap\Z$.

\subsection{Independent perturbations}
In this subsection we assume that the variables ${\Mpert}_{\dxp}$ are independent and identically distributed. 
\begin{theorem}
\label{thm:squaretrans}
The square modulus of the transmission coefficient $\tau = |T|^2$ behaves as a diffusion process with the infinitesimal generator 
\begin{equation}\begin{split}
{\mathcal L} =\gamma \big[ \tau^2(1-\tau) \partial_\tau^2 - \tau^2 \partial_\tau\big], 
\label{eq:defL}
\end{split}\end{equation}
starting from $\tau_0=0$,
where 
\begin{equation}\begin{split}
\label{def:gamma:ind2}
\gamma(\idf) = \frac{{\idf}^4 \sigma^2}{4 \sin^2 k(\idf) } ,
\end{split}\end{equation}
and $k(\idf)$ is defined by (\ref{eq:defkomegas}).
\end{theorem}

This theorem is proved in Appendix \ref{sec:proofthm:squaretrans}.
The form of the infinitesimal generator is similar to the one obtained for the square modulus of the transmission coefficient  of the one-dimensional wave equation in random medium \cite[Theorem  7.3]{book}, except for the frequency dependence which is here original and which follows from the particular dispersion relation of the discrete lattice.

Using the results of \cite[Section 7.1.5]{book} we obtain that, for any $n\geq 1$:
\begin{equation}\begin{split}
\label{eq:momT}
\mathbb{E} \big[ |T|^{2n} \big] 
&=e^{-\frac{\gamma {\Ld}}{4}}
\int_0^\infty e^{- \gamma {\Ld} s^2} \frac{2 \pi s \sinh(\pi s)}{\cosh^2(\pi s)} \phi_n(s) ds ,
\end{split}\end{equation}
with
\begin{equation}\begin{split}
\phi_1(s)=1,\qquad \phi_n(s)= \prod_{j=1}^{n-1} \frac{s^2 +(j-\frac{1}{2})^2}{j^2}, \quad n \geq 2.
\label{eq:defphins}
\end{split}\end{equation}
As a result of \eqref{eq:momT}, we can obtain the mean transmission $\mathbb{E} \big[ |T|^{2} \big] $
and its variance ${\rm Var}\big( |T|^{2} \big)= \mathbb{E} \big[ |T|^{4} \big] -\mathbb{E} \big[ |T|^{2} \big]^2$.

\begin{figure}[ht!]
(a)\includegraphics[width=.45\textwidth]{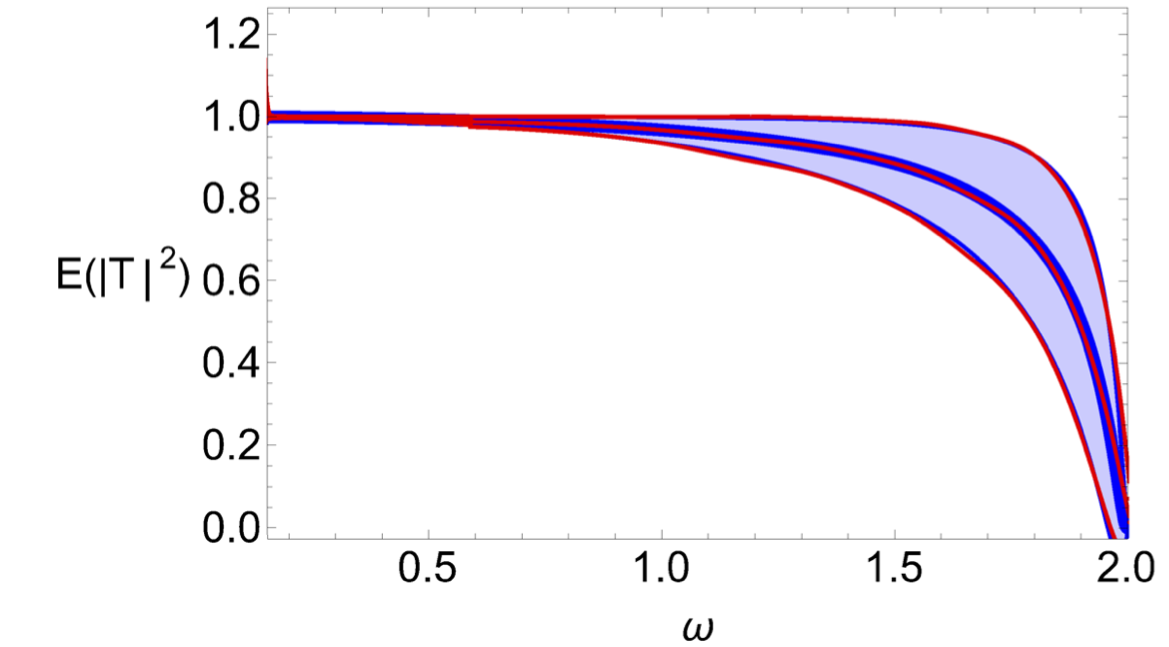}(b)\includegraphics[width=.45\textwidth]{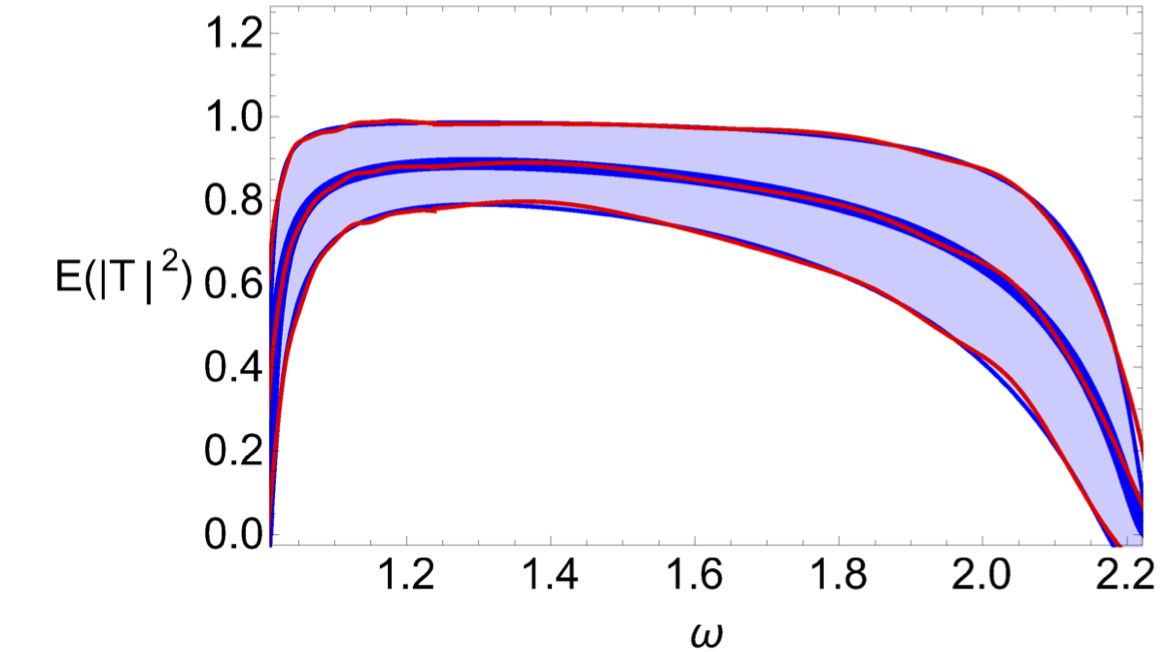}
\caption[ ]{Mean curve $\mathbb{E}[|T|^2]$ (thick curve) with spread $\pm {\rm Std} (|T|^2)$ vs frequency $\fq$.
(a) $\Ks=0.02$ (b) $\Ks=1.01$.
Blue: asymptotic formula \eqref{eq:momT}, Red: emprirical averages of numerical simulations using the exact expression \eqref{Tfullmatched} described in Appendix \ref{Tfullmatchedapp}. Ensemble details are the same ones as in Fig. \ref{freqdomainonedimlatticerandommasspertb}. In the absence of mass perturbation on $[1,\Ld]\cap\Z$, we have $|T|^2\equiv1.$}
\label{onedimlatticerandom12}
\end{figure}

In Figure \ref{onedimlatticerandom12}, we compare the empirical averages of numerical simulations with the theoretical predictions for the expectation $\mathbb{E}[|T|^2]$ and the standard deviation ${\rm Std}( |T|^{2} )={\rm Var}( |T|^{2} )^{1/2}$. The numerical simulations are based on the exact solution \eqref{Tfullmatched} and the theoretical predictions are based on \eqref{eq:momT}.
We obtain excellent agreement which confirms the accuracy of the asymptotic approach.
%The blue dots in Fig.~\ref{Cases} correspond to \eqref{nonmatchedcase1} (both plots differ by sign of ${\Mpt}_0$).
We can observe that the behavior of the transmittance close to the left endpoint of the propagative band $(\sqrt{\Ks},\sqrt{\Ks+4})$ is very different in the cases $\Ks=0$ and $\Ks\neq 0$.
The transmittance goes to zero as $\omega \to \sqrt{\Ks}$ when $\Ks>0$ and it goes to one when $\Ks=0$.

\subsection{Correlated perturbations}
The previous results can be extended to the case where the variables ${\Mpert}_{\dxp}$ are identically distributed with mean zero, variance $\sigma^2$, and integrable covariance function:
\begin{equation}\begin{split}
\mathbb{E}[{\Mpert}_{\dxp} {\Mpert}_{{\dxp}'}]=\sigma^2 \Gamma({\dxp}-{\dxp}'),\quad \dxp, \dxp'\in\Z,
\label{eq:covm}
\end{split}\end{equation}
with $\sigma \ll 1$, and ${\Ld}$ is of the order of $\sigma^{-2}$.
The function $\Gamma$ is assumed to be integrable $\sum_{j \in \mathbb{Z}} |\Gamma(j)|<+\infty$. We can then apply the diffusion approximation theory set forth in \cite[Chapter 6]{book} 
(other approaches based on Duhamel series expansions exist but will not be used here \cite{balgu15}).
We find that $\tau = |T|^2$ behaves as a diffusion process with the infinitesimal generator (\ref{eq:defL})
where 
\begin{equation}\begin{split}
\label{eq:defgamma2}
\gamma(\idf)= \frac{{\idf}^4 \sigma^2}{4 \sin^2 k(\idf) } \check{\Gamma}(2k(\idf)), \qquad \check{\Gamma}(k) = \Gamma(0)+2\sum_{j=1}^\infty \cos(kj) \Gamma(j).
\end{split}\end{equation}
The moments of $|T|^2$ are still given by (\ref{eq:momT}) with the new expression (\ref{eq:defgamma2}) of $\gamma$.
Note that $\check{\Gamma}(2k)$ is non negative by Wiener-Khintchine theorem but it may be non-monotoneous as a function of $k$. In other words the spatial correlation of the variables ${\Mpert}_{\dxp}$ has a non-trivial impact onto the parameter $\gamma(\idf)$.
As a consequence, the observation of $\gamma(\idf)$ for $\idf \in (\sqrt{\Ks},\sqrt{\Ks+4})$ makes it possible to retrieve $\check{\Gamma}(k)$ for $k\in (0,2\pi)$ which characterizes the correlation function $\Gamma({\dxp})$. This opens the way towards an original method to estimate the statistics of random mass perturbations by measuring the frequency-dependent transmission coefficient of a section of a discrete lattice.

\subsection{Scattering in non-matched medium} 
The previous results can be revisited in the case where
\begin{align}
{\Mpert}_{\dxp} = {\Mpt}_0 & \mbox{ for } {\dxp} \leq 0 ,\\
{\Mpert}_{\dxp} = {\Mpt}_1 & \mbox{ for } {\dxp} > {\Ld},
\end{align}
and $\dxp\in\Z.$
This means that the masses in the two unperturbed half-spaces may be different from the average mass in the perturbed section $[1,{\Ld}]\cap\Z$. 
We may then anticipate that the boundaries of the perturbed section can generate themselves reflections.

We introduce the wavenumbers ${\idk}_0$ and ${\idk}_1$ solutions of the dispersion relations
\begin{equation}\begin{split}
\label{eq:dispersionrelation:j}
-(1+{\Mpt}_j) {\idf}^2 =2\cos {\idk}_j -2-\Ks , \quad j=0,1.
\end{split}\end{equation}
Here we assume the regime is propagative, i.e. the frequency ${\idf}$ is such that $(2+\Ks-{\idf}^2(1+{\Mpt}_j))/(2) \in (-1,1)$ for $j=0,1$ so that there is a unique solution ${\idk}_j\in (0,\pi)$ to (\ref{eq:dispersionrelation:j}).
If a right-going input wave is incoming from the left, the time-harmonic field in the two unperturbed half-spaces has the form
\begin{align}
\hat{u}_{\dxp} &= e^{i{\idk}_0{\dxp}}+R e^{-i{\idk}_0{\dxp}} \quad \mbox{ for } {\dxp} \leq 0,\\
\hat{u}_{\dxp} &= T e^{i{\idk}_1{\dxp}}\quad  \mbox{ for } {\dxp} >{\Ld},
\label{eq:formuyafetrL}
\end{align}
and $\hat{u}$ satisfies (\ref{eq:th}) for $ 0 < {\dxp} \leq {\Ld}$ with $\dxp\in\Z$.

We assume that the random variables ${\Mpert}_{\dxp}$ are identically distributed with mean zero, variance $\sigma^2$, and integrable covariance function (\ref{eq:covm}).
We can give explicit formulas in two special cases.

1) If ${\Mpt}_1=0$, then we get
\begin{equation}\begin{split}
\mathbb{E}[|T|^2 ] &=
\frac{2\sin^2 {\idk}_0}{1-\cos(k+{\idk}_0)}
\sum_{m=0}^\infty \Big( \frac{1-\cos(k-{\idk}_0)}{1-\cos(k+{\idk}_0)}\Big)^m \Big( \mathbb{E} [|\tilde{R}|^{2m}] - \mathbb{E} [|\tilde{R}|^{2m+2}] \Big)
\\ &= \frac{2\sin^2 {\idk}_0}{1-\cos(k+{\idk}_0)}
\sum_{n=0}^\infty (-1)^n \mathbb{E} [|\tilde{T}|^{2n+2}] \Big[ \sum_{m=n}^\infty \binom{m}{n} \Big( \frac{1-\cos(k-{\idk}_0)}{1-\cos(k+{\idk}_0)}\Big)^m \Big] ,
\label{nonmatchedcase1}
\end{split}\end{equation}
where $\mathbb{E}[|\tilde{T}|^{2n}]$ is given by (\ref{eq:momT}).
When there is no random mass perturbation,
we have simply $|T|^2 = \frac{2 \sin^2 k_0}{1 - \cos(k+k_0)}$ (see also \eqref{nonmatchexactTperf}).
We have similarly
\begin{equation}\begin{split}
\mathbb{E}[|T|^4 ]
=& \Big( \frac{2\sin^2 {\idk}_0}{1-\cos(k+{\idk}_0)} \Big)^2
\sum_{n=0}^\infty (-1)^n \mathbb{E} [|\tilde{T}|^{2n+4}]\\
&\times \Big[ \sum_{m=n}^\infty \binom{m}{n} (1+m)^2 \Big( \frac{1-\cos(k-{\idk}_0)}{1-\cos(k+{\idk}_0)}\Big)^m \Big] ,
\label{nonmatchedcase1a}
\end{split}\end{equation}
which makes it possible to get ${\rm Var}(|T|^2) = \mathbb{E}[|T|^4]- \mathbb{E}[|T|^2]^2$.

2) If ${\Mpt}_0=0$, then
we have
\begin{equation}\begin{split}
\mathbb{E}[|T|^2] =\frac{\sin k}{\sin {\idk}_1} 
e^{-\frac{\gamma {\Ld}}{4}}
\int_0^\infty e^{- \gamma {\Ld} s^2} \frac{2 \pi s \sinh(\pi s)}{\cosh^2(\pi s)} P_{-1/2+is}\Big(
\frac{1-\cos k \cos {\idk}_1}{\sin k \sin {\idk}_1}\Big) ds .
\label{nonmatchedcase2}
\end{split}\end{equation}
When there is no random mass perturbation,
we have simply $|T|^2 = \frac{2 \sin^2 k}{1 - \cos(k+k_1)}$ (see also \eqref{nonmatchexactTperf}).
More generally, for any $n\geq 1$,\begin{equation}\begin{split}
\mathbb{E}[|T|^{2n}] =\Big( \frac{\sin k}{\sin {\idk}_1} \Big)^n
e^{-\frac{\gamma {\Ld}}{4}}
\int_0^\infty e^{- \gamma {\Ld} s^2} \frac{2 \pi s \sinh(\pi s)}{\cosh^2(\pi s)} \phi_n(s)\\
\times
P_{-1/2+is}\Big(
\frac{1-\cos k \cos {\idk}_1}{\sin k \sin {\idk}_1}\Big) ds ,
\label{nonmatchedcase2nn}
\end{split}\end{equation}
where $\phi_n(s)$ is defined by (\ref{eq:defphins})
and
\begin{equation}\begin{split}
P_{-1/2+is}(\eta) = \frac{\sqrt 2}{\pi} \cosh(\pi s) \int_{0}^{\infty}
\frac{\cos(s \tau)}{\sqrt{\cosh (t + \eta)}} \, d t.
\label{Phalfeta}
\end{split}\end{equation}

\begin{remark}
It is possible to address the general case where both ${\Mpt}_0$ and ${\Mpt}_1$ are not zero, 
but the expressions become complicated.
\end{remark}

\begin{remark}
It is possible to address the case where ${\idf}$ is outside the pass band on the right half-space, i.e. $(2+\Ks-{\idf}^2(1+{\Mpt}_1))/(2) \not\in (-1,1)$. 
The results are given in Subsection \ref{sec:proofnonmatched}. As can be expected the reflection coefficients has modulus one since the wave cannot propagate in the right half-space.
\end{remark}

\begin{figure}[htb!]
\centering
(a)\includegraphics[width=.45\textwidth]{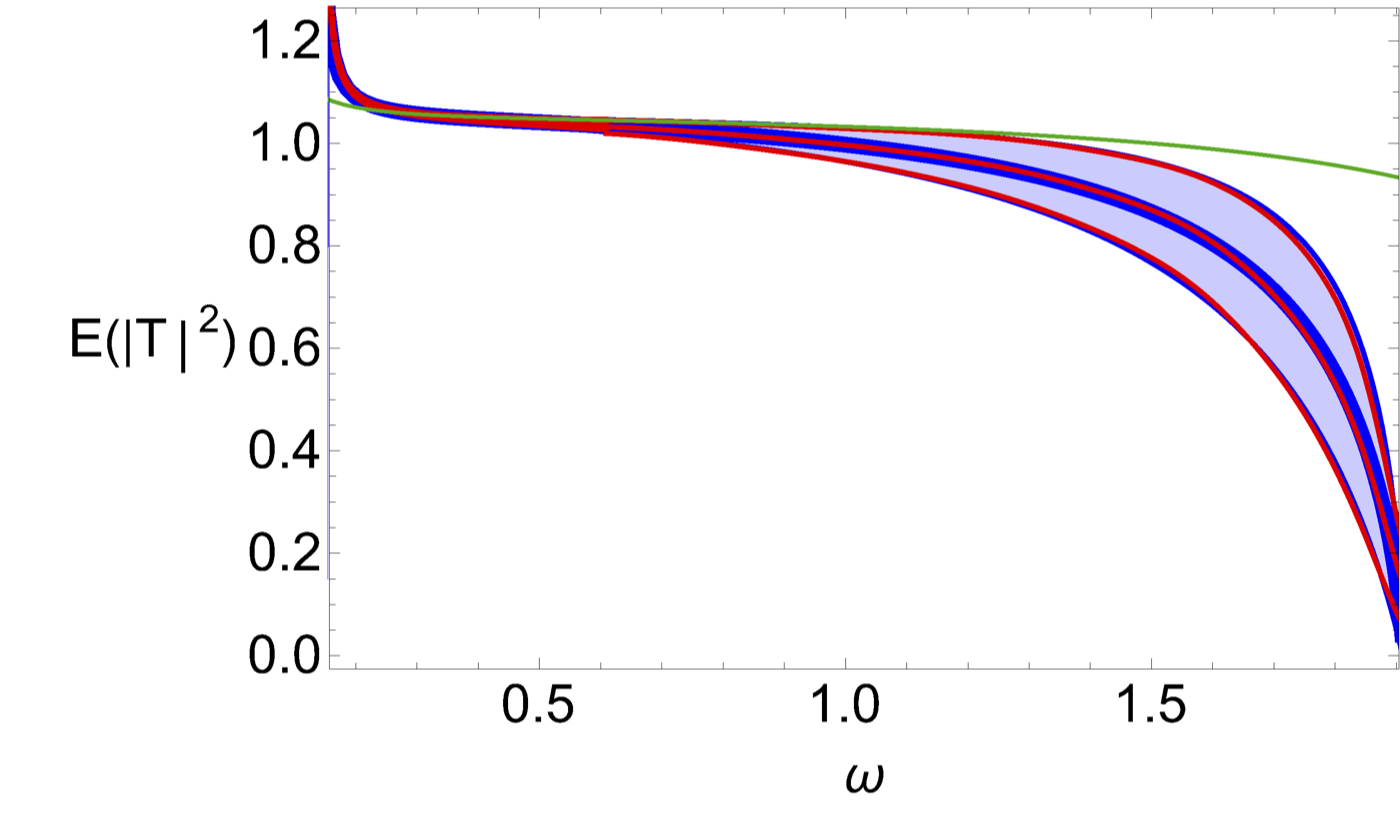}(b) 
\includegraphics[width=.45\textwidth]{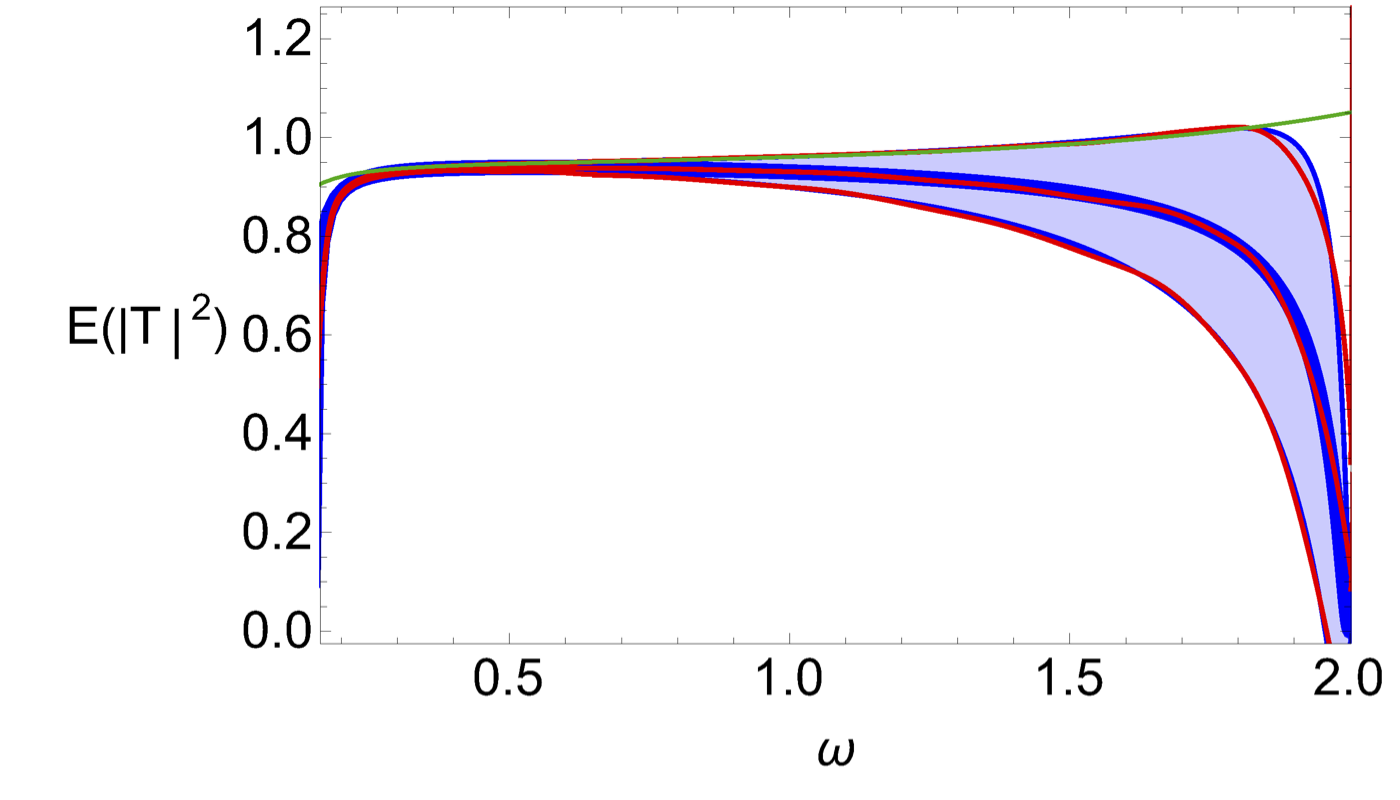}\\
(c)\includegraphics[width=.45\textwidth]{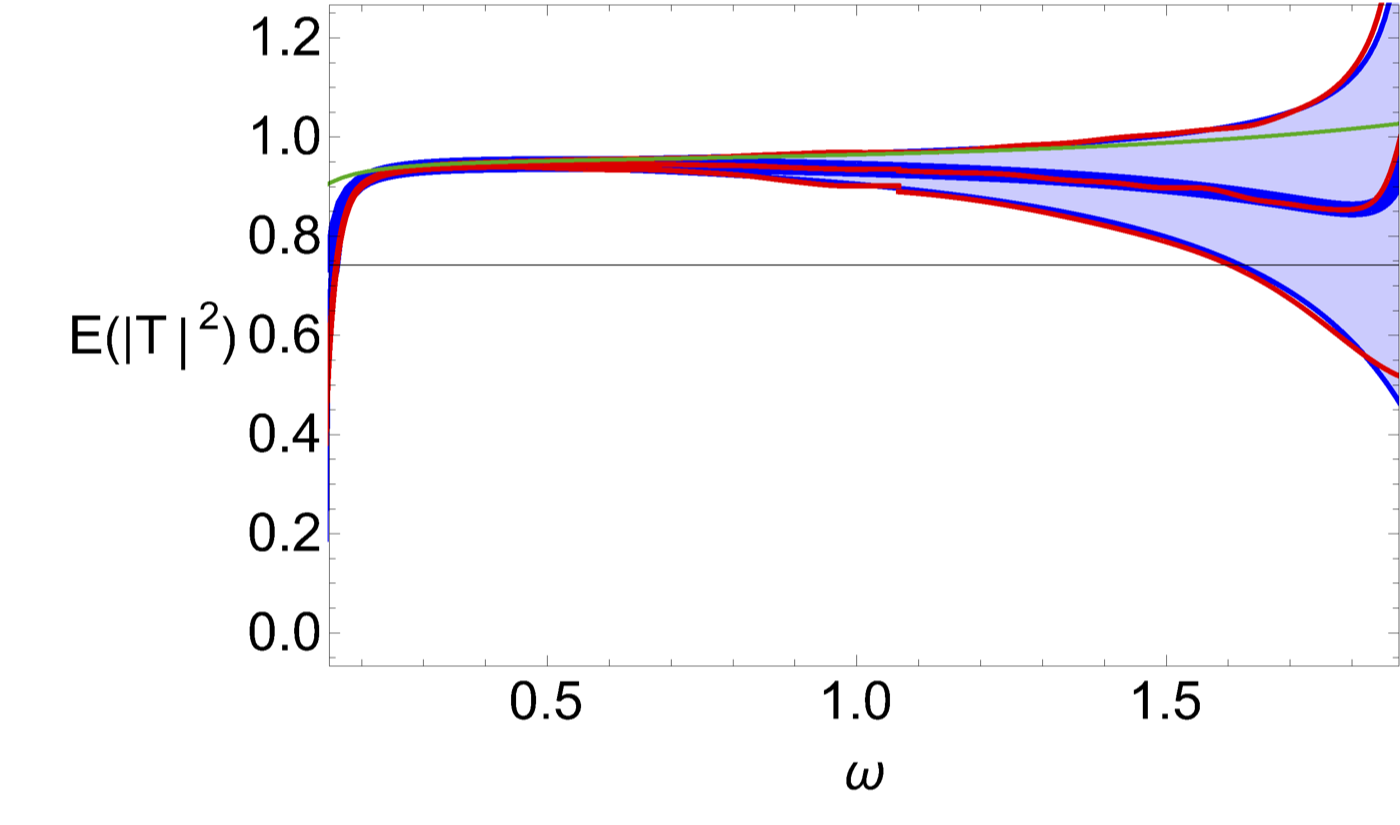}(d)
\includegraphics[width=.45\textwidth]{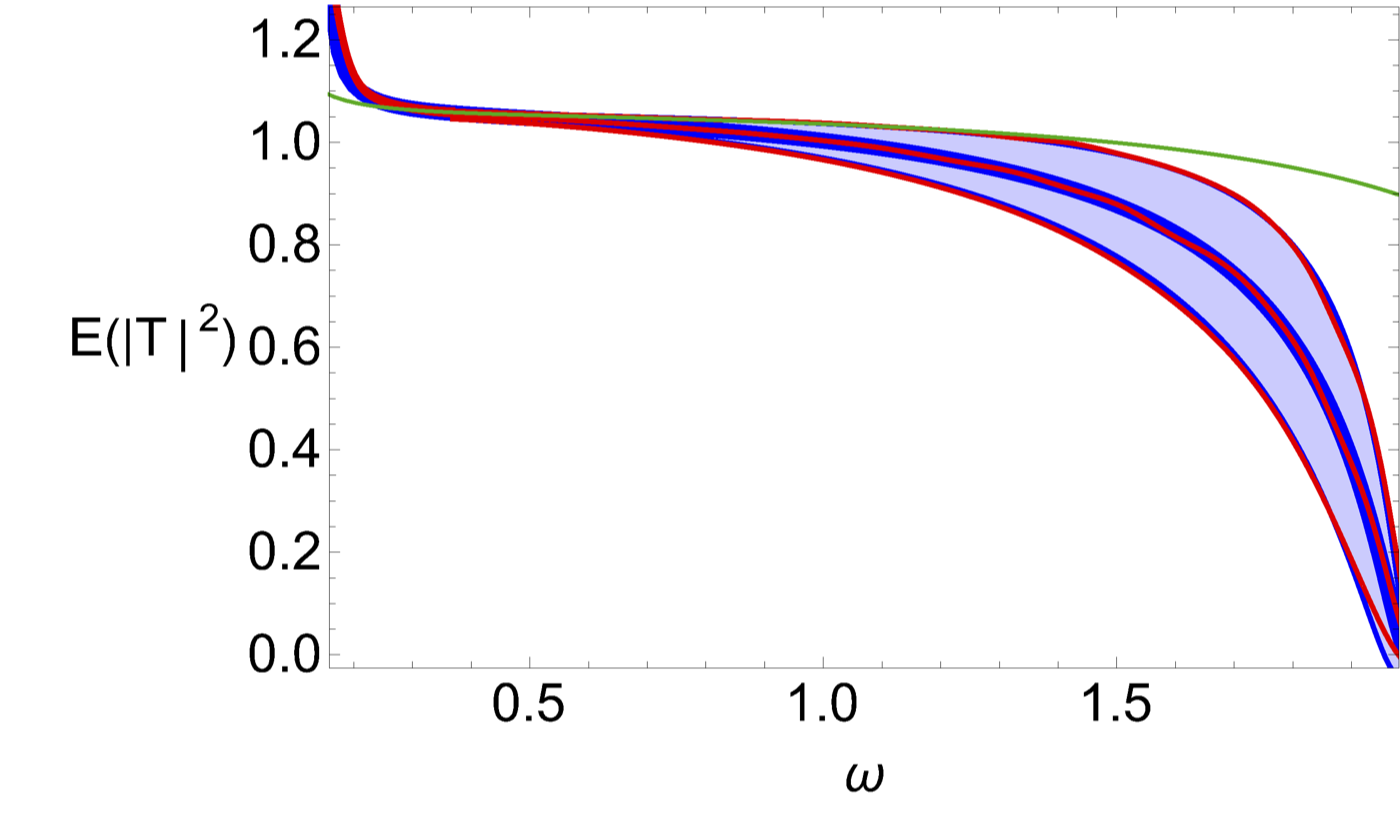}
\caption{Mean curve $\mathbb{E}[|T|^2]$ (thick curve) with spread $\pm {\rm Std} (|T|^2)$ vs frequency $\fq$.
(a) $\Delta_0=0.1, \Delta_1=0$,
(b) $\Delta_0=-0.1, \Delta_1=0$,
(c) $\Delta_0=0, \Delta_1=0.1$,
(d) $\Delta_0=0, \Delta_1=-0.1$.
In all cases ${\Nd}=40, \sigma=0.05, \Ks=0.02.$
Blue: asymptotic formulas. Red: empirical averages based on the exact expression  \eqref{nonmatchexactT} for (a), (b) and \eqref{nonmatchexactT2a} for (c), (d); the ensemble size is $151$. Green: formula of $|T|^2$ \eqref{nonmatchexactTperf} for the case \eqref{nonmatchmass}  without mass perturbation on $[1,\Ld]\cap\Z$;
similarity between (a), (c) (and (b), (d)) is due to the symmetry in the expression of $|T|^2$ with respect to index $0$ and $1$ in the absence of mass perturbation. }
\label{Cases}\end{figure}

In Figure \ref{Cases} we compare the empirical averages of numerical simulations with the theoretical predictions for the expectation $\mathbb{E}[|T|^2]$ and the standard deviation ${\rm Std}(|T|^2)$ for different mismatched media. The numerical simulations are based on the exact solution \eqref{nonmatchexactT} for $\Delta_0\neq 0$, $\Delta_1= 0$ and \eqref{nonmatchexactT2a} for $\Delta_0=0$, $\Delta_1\neq 0$. The theoretical predictions are based on \eqref{nonmatchedcase1} and \eqref{nonmatchedcase1a} for $\Delta_0\neq 0$, $\Delta_1= 0$ and 
on \eqref{nonmatchedcase2nn} for $\Delta_0=0$, $\Delta_1\neq 0$. 
We obtain excellent agreement which confirms the accuracy of the asymptotic approach.
We can observe that the medium mismatch has a dramatic impact on the transmittance, in particular close to the endpoints of the propagative band $(\sqrt{\Ks},\sqrt{\Ks+4})$.

%\textcolor{red}{can we comment on the large variation in the upper-centre of band vs thinning near the band-edge in (a), (b), (d) vs (c) of Fig. \ref{Cases}? Add a remark on the formula for $|T|^2$ for just mass mismatch with out perturbation, add that curve in figure 6. See Remark \ref{remarkTmisperf}.}

\section{Proofs of Preliminary results}
\label{sec:proofprel}
\begin{proof}[Proof of Lemma \ref{lemma1}]
Since $k(\fq)$ satisfies $-\fq^2 = 2 (\cos(k(\fq))-1)$, we find that $\hat{\su}_{\dxp}(\fq)$ satisfies 
\begin{equation}\begin{split}
\hat{\su}_{{\dxp}+1}(\fq)+\hat{\su}_{{\dxp}-1}(\fq)-2\hat{\su}_{\dxp}(\fq)
= 2 (\cos(k(\fq))-1) \hat{\su}_{\dxp}(\fq)
=-\fq^2 \hat{\su}_{\dxp}(\fq)
\end{split}\end{equation} 
for any $\dxp \in \Z$, hence ${\su}_{\dxp}(t)$
satisfies (\ref{TDeq}).
It remains to show that ${\su}_{{\dxp}}(t)$ satisfies the appropriate initial conditions. 
On the one hand (using the change of variable $\fq \mapsto k(\fq)$) we have
\begin{equation}\begin{split}
{\su}_{{\dxp}}(0) &= \frac{1}{2\pi} \int_0^2 \hat{\su}_{{\dxp}}(\fq)  d\fq + c.c. \\
&= 
\frac{2}{\pi} \int_0^2 \frac{\cos(k(\fq)(\dxp-\dxp_0))}{\sqrt{4-\fq^2}}  d\fq =
\frac{1}{\pi} \int_{0}^\pi \cos(k (\dxp-\dxp_0))dk  =\delta_{\dxp , \dxp_0}.
\end{split}\end{equation}
On the other hand we have
\begin{equation}\begin{split}
\dot{\su}_{{\dxp}}(0) = - \frac{i}{2\pi} \int_0^2 \fq \hat{\su}_{{\dxp}}(\fq)  d\fq + c.c. = 0 ,
\end{split}\end{equation}
which completes the proof.
\end{proof}

\begin{proof}[Proof of Lemma \ref{lemma2}]
For $\dxp \gg 1, \dxp\in\Z,$, the expression 
\begin{equation}\begin{split}
{\su}_{\dxp_0+\dxp}(\alpha \dxp) =& \frac{1}{2\pi} \int_0^2 \frac{1}{\sqrt{4-\fq^2}}
e^{i ( k(\fq) \dxp- \fq \alpha \dxp)} d\fq\\
& +  \frac{1}{2\pi} \int_0^2 \frac{1}{\sqrt{4-\fq^2}}
e^{i ( - k(\fq) \dxp - \fq \alpha \dxp) } d\fq + 
c.c.
\end{split}\end{equation}
involves the value of an integral (in $\fq$) with a fast phase $[\pm k(\fq) - \fq \alpha] \dxp$ 
and it can be evaluated by the stationary phase method. 
We compute
\begin{equation}\begin{split}
k'(\fq) = \frac{2}{\sqrt{4-\fq^2}} ,\qquad
k''(\fq) = \frac{2\fq}{ (4-\fq^2)^{3/2}} ,
\end{split}\end{equation}
The phase $ [ k(\fq)  - \fq \alpha]\dxp$  has a unique stationary point at $\fq_\alpha \in (0,2) $ such that $k'(\fq_\alpha ) =\alpha$, i.e.
\begin{equation}
    \label{def:omegaalpha2}
\fq_\alpha =  \frac{2\sqrt{\alpha^2-1}}{\alpha} , 
\end{equation}
for which we have $k''(\fq_\alpha) =  {\alpha^2 \sqrt{\alpha^2-1}}/{2}$.
We then obtain (after the change of variable $\fq=\fq_\alpha + \dxp^{-1/2} s $)
\begin{align}
{\su}_{\dxp_0+\dxp}(\alpha \dxp)  &= \frac{1}{2\pi \sqrt{4-\fq_\alpha^2} \sqrt{\dxp}} 
e^{i [k(\fq_\alpha)  -\fq_\alpha \alpha] \dxp} 
 \int_{-\infty}^\infty e^{i \frac{k''(\fq_\alpha) }{2} s^2} ds 
 +c.c. +o(\frac{1}{\sqrt{\dxp}}) \notag \\
&= \frac{1}{2\pi \sqrt{4-\fq_\alpha^2} \sqrt{\dxp}} 
e^{i [k(\fq_\alpha)  -\fq_\alpha \alpha]\dxp } 
\frac{\sqrt{2\pi} e^{i \pi/4} }{\sqrt{k''(\fq_\alpha)}} + c.c. +o(\frac{1}{\sqrt{\dxp}}),
\end{align}
which gives (\ref{eq:asyhomouy1}).
\end{proof}

\begin{proof}[Proof of Lemma \ref{lemma3}]
We now look for an asymptotic expansion of field ${\su}_{\dxp_0+\dxp} (t)$ around time $\dxp$ with $\dxp>0, \dxp\in\Z$:
\begin{equation}\begin{split}
 {\su}_{\dxp_0+\dxp} ( \dxp+\sqrt[3]{\dxp} \beta) 
& = \frac{1}{2\pi} \int_0^2 \frac{1}{\sqrt{4-\fq^2}}
e^{i ( k(\fq) \dxp- \fq  \dxp - \fq \beta \sqrt[3]{\dxp})} d\fq \\
&+  \frac{1}{2\pi} \int_0^2 \frac{1}{\sqrt{4-\fq^2}}
e^{i ( - k(\fq) \dxp -\fq \dxp - \fq \beta \sqrt[3]{\dxp}) } d\fq + 
c.c. 
 \end{split}\end{equation}
For $\dxp \gg 1$, we find that the phase $ [k(\fq) \dxp -\fq  \dxp - \fq \beta \sqrt[3]{\dxp}]$  has a unique stationary point at $\fq_0=0 $ such that $k'(\fq_0)=1$ and it is localized at the border of the interval $(0,2)$. 
We then obtain (using $k'''(0)=1/4$ and the change of variable $\fq=\dxp^{-1/3} s$)
\begin{align}
 {\su}_{\dxp_0+\dxp} ( \dxp+\sqrt[3]{\dxp} \beta) &= \frac{1}{4\pi \sqrt[3]{\dxp}}  
 \int_{0}^\infty e^{i \frac{ k'''(0) }{6} s^3- i s \beta } ds 
 +c.c. +o(\frac{1}{\sqrt[3]{\dxp}}) \notag \\
&= \frac{1}{\pi \sqrt[3]{\dxp}} \int_0^\infty \cos \big(  \frac{s^3}{3} - 2 \beta s\big) ds +o(\frac{1}{\sqrt[3]{\dxp}}) ,
\end{align}
which completes the proof since $Ai({\dxp}) = \frac{1}{\pi}\int_0^\infty \cos \big( {\dxp} s +\frac{s^3}{3}\big) ds$.
\end{proof}

\begin{proof}[Proof of Lemma \ref{lemma4}]
For $\dxp \gg 1$, the expression 
\begin{align}
\nonumber
{\su}_{\dxp_0+\dxp}(\alpha \dxp) =& \frac{1}{4\pi} \int_0^\infty \hat{c}(\fq)
e^{i ( k(\fq) \dxp- \fq \alpha \dxp)} d\fq \\
&+  \frac{1}{4\pi} \int_0^\infty \hat{c}(\fq)
e^{i ( - k(\fq) \dxp - \fq \alpha \dxp) } d\fq + 
c.c.
\end{align}
is the value of an integral (in $\fq$) with a fast phase $[\pm k(\fq) - \fq \alpha] \dxp$ 
and it can be evaluated by the stationary phase method. 
We compute
\begin{equation}\begin{split}
k'(\fq) = \frac{2 \fq}{\sqrt{\fq^2-{\Ks}}\sqrt{4+{\Ks} - \fq^2 }} ,\qquad
k''(\fq) = \frac{2(\fq^4 - \fq_s^4)}{(\fq^2-{\Ks})^{3/2} (4+{\Ks}  - \fq^2)^{3/2}} .
\end{split}\end{equation}
If $\alpha < \alpha_s$ then the phase $ [ k(\fq)  - \fq \alpha]\dxp$ does not have any stationary point.
For any $\alpha \in (\alpha_s,+\infty)$
the phase $ [ k(\fq)  - \fq \alpha]\dxp$  has two stationary points at $\fq_{\alpha}^{\pm} \in (\sqrt{{\Ks}},\sqrt{4+{\Ks}}) $ such that $k'(\fq_{\alpha}^{\pm} ) =\alpha$.
We have in fact $\fq_{\alpha}^{+} \in (\fq_s,\sqrt{4+{\Ks}})$ and $\fq_{\alpha}^{-} \in (\sqrt{{\Ks}},\fq_s)$,
and
\begin{equation}
\begin{split}
k''(\fq_{\alpha}^{\pm}) = \frac{\alpha^3 ( {\fq_{\alpha}^{\pm}}^4-\fq_s^4)}{4{\fq_{\pm ,\alpha}}^3} , 
\end{split}
\end{equation}
which is positive for $\fq_{\alpha}^{+}$ and negative for $\fq_{\alpha}^{-}$.
We then obtain
\begin{align}
{\su}_{\dxp_0+\dxp}(\alpha \dxp)  &= \frac{1}{4\pi \sqrt{\dxp} }\hat{c}(\fq_{\alpha}^{+}) 
e^{i [k(\fq_{\alpha}^{+})  -\fq_{\alpha}^{+} \alpha] \dxp} 
 \int_{-\infty}^\infty e^{i \frac{k''(\fq_{\alpha}^{+}) }{2} s^2} ds \notag\\
 &\quad 
 +\frac{1}{4\pi \sqrt{\dxp} }\hat{c}(\fq_{\alpha}^{-}) 
e^{i [k(\fq_{\alpha}^{-})  -\fq_{\alpha}^{-} \alpha] \dxp} 
 \int_{-\infty}^\infty e^{i \frac{k''(\fq_{\alpha}^{-}) }{2} s^2} ds 
 +c.c. +o(\frac{1}{\sqrt{\dxp}})\notag  \\
&= \frac{ \hat{c}(\fq_{\alpha}^{+}) }{\sqrt{2\pi k''(\fq_{\alpha}^{+}) \dxp}}
\cos\Big( \frac{\pi}{4} +  [k(\fq_{\alpha}^{+})  -\fq_{\alpha}^{+} \alpha] \dxp\Big)
\notag\\
&\quad + \frac{  \hat{c}(\fq_{\alpha}^{-}) }{2\sqrt{\pi |k''(\fq_{\alpha}^{-})| \dxp}}
\cos\Big( -\frac{\pi}{4} +  [k(\fq_{\alpha}^{-})  -\fq_{\alpha}^{-} \alpha] \dxp\Big)
+o(\frac{1}{\sqrt{\dxp}}),
\end{align}
which gives (\ref{eq:asyhomouy1:ks}) using $\hat{c}(\fq_{\alpha}^{\pm}) =  \alpha$.
\end{proof}

\begin{proof}[Proof of Lemma \ref{lemma5}]
We now look for an asymptotic expansion of field ${\su}_{\dxp_0+\dxp} (t)$ around time $\alpha_s \dxp$ with $\dxp>0, \dxp\in\Z$:
\begin{equation}\begin{split}
 {\su}_{\dxp_0+\dxp} (\alpha_s \dxp+ \beta \sqrt[3]{\dxp}) 
 =& \frac{1}{4\pi} \int_0^\infty \hat{c}(\fq)
e^{i ( k(\fq) \dxp- \fq \alpha_s \dxp - \fq \beta \sqrt[3]{\dxp})} d\fq\\
& +  \frac{1}{4\pi} \int_0^\infty\hat{c}(\fq)
e^{i ( - k(\fq) \dxp - \fq \alpha_s \dxp - \fq \beta \sqrt[3]{\dxp}) } d\fq + 
c.c. 
 \end{split}\end{equation}
For $\dxp \gg 1$, we find that the phase $ [k(\fq) \dxp -\fq \alpha_s \dxp - \fq \beta \sqrt[3]{\dxp}]$  has a unique stationary point at $\fq_s $ such that $k'(\fq_s)=\alpha_s$. 
It satisfies $k''(\fq_s)=0$ and $\fq_s$ is localized in the interior of the interval $(\sqrt{{\Ks}},\sqrt{4+{\Ks}})$. 
We then obtain by the change of variable $\fq=\fq_s+\dxp^{-1/3} s$ (and using $k'''(\fq_s)=\alpha_s^3$, $\hat{c}(\fq_s)=\alpha_s$)
\begin{align}
 {\su}_{\dxp_0+\dxp} ( \alpha_s \dxp+\beta \sqrt[3]{\dxp}) =&  \frac{\hat{c}(\fq_s)}{4\pi \sqrt[3]{\dxp} } e^{i [k(\fq_s) - \fq_s \alpha_s] \dxp - i \fq_s \beta \sqrt[3]{\dxp}}
 \int_{-\infty}^\infty e^{i \frac{ k'''(\fq_s)  }{6} s^3- i s \alpha_s } ds \notag\\
& +c.c. +o(\frac{1}{\sqrt[3]{\dxp}}) \notag \\
=& \frac{\sqrt[3]{2}}{2\pi  \sqrt[3]{\dxp}} \int_0^\infty \cos \big(  \frac{s^3}{3} - \sqrt[3]{2} \frac{ \beta}{\alpha_s} s\big) ds e^{i [k(\fq_s) - \fq_s \alpha_s] \dxp - i \fq_s \beta \sqrt[3]{\dxp}}\notag\\
&+c.c. + o(\frac{1}{\sqrt[3]{\dxp}}) ,
\end{align}
which completes the proof of (\ref{eq:asyhomouy2:ks}).
\end{proof}

\section{Proof of Theorem \ref{thm:squaretrans}}
\subsection{Scattering in matched medium \texorpdfstring{${\Mpt}_0={\Mpt}_1=0$}{}}
\label{sec:proofthm:squaretrans}
When $\Mpert_{\dxp} \equiv 0$, the solution is of the form
\begin{equation}
\hat{u}_{\dxp} = \alpha e^{i k {\dxp}}+\beta e^{-ik{\dxp}} ,\quad \dxp\in\Z,
\end{equation}
with $k$ solution of the dispersion relation
\begin{equation}
\label{eq:dispersionrelation}
-\idf^2 =2\cos k -2 - \Ks.
\end{equation}
Here we assume the regime is propagative, i.e. the frequency $\idf$ is such that $(2+\Ks-\idf^2)/(2) \in (-1,1)$ so that there is a unique solution $k\in (0,\pi)$ to (\ref{eq:dispersionrelation}).

When $\Mpert_{\dxp} =0$ for ${\dxp} \leq 0$ and for ${\dxp}>{\Ld}$
and a right-going input wave is incoming from the left, the solution has the form
\begin{align}
\hat{u}_{\dxp} &= e^{ik{\dxp}}+R e^{-ik{\dxp}} \quad \mbox{ for } {\dxp} \leq 0,\\
\hat{u}_{\dxp} &= T e^{ik{\dxp}}\quad  \mbox{ for } {\dxp} >{\Ld},
\end{align}
and $\hat{u}$ satisfies (\ref{eq:th}) for $ 0 < {\dxp} \leq {\Ld}$ with $\dxp\in\Z$.

We introduce 
\begin{align}
\alpha_{\dxp} &= \frac{e^{-ik{\dxp}}}{2i\sin k} \big( \hat{u}_{{\dxp}+1}-e^{-ik} \hat{u}_{\dxp} \big) ,\\ 
\beta_{\dxp} &= - \frac{e^{ik{\dxp}}}{2i\sin k} \big( \hat{u}_{{\dxp}+1}-e^{ik} \hat{u}_{\dxp} \big) .
\end{align}
We then have $(\alpha_{\dxp},\beta_{\dxp})=(1,R)$ for ${\dxp} \leq 0$, $(\alpha_{\dxp},\beta_{\dxp})=(T,0)$ for ${\dxp}>{\Ld}$, and 
\begin{equation}
\hat{u}_{\dxp} = \alpha_{\dxp} e^{ik{\dxp}}+\beta_{\dxp} e^{-ik{\dxp}} ,
\end{equation}
for any ${\dxp}\in \mathbb{Z}$.

After some algebra, using the fact that
\begin{equation}
\hat{u}_{{\dxp}+1}-\hat{u}_{\dxp} = \alpha_{\dxp} e^{ik{\dxp}}(e^{ik}-1)+\beta_{\dxp} e^{-ik{\dxp}}(e^{-ik}-1), 
\end{equation}
we find that $(\alpha_{\dxp},\beta_{\dxp})$ satisfies the system:
\begin{align}
\alpha_{{\dxp}-1} &=\alpha_{\dxp} -\frac{i \idf^2}{2 \sin k} \Mpert_{\dxp} \big(\alpha_{\dxp} +\beta_{\dxp} e^{-2ik{\dxp}}\big),\\
\beta_{{\dxp}-1} &=\beta_{\dxp} +\frac{i \idf^2}{2 \sin k} \Mpert_{\dxp} \big(\alpha_{\dxp} e^{2ik{\dxp}}+\beta_{\dxp} \big),
\end{align}
for $0 < {\dxp} \leq {\Ld}$, with the boundary conditions
\begin{align}
\alpha_{\Ld}=T,  \quad \beta_{\Ld}=0, \quad \alpha_0=1, \quad \beta_0=R.
\end{align}

\subsubsection{Expressions of the reflection and transmission coefficients}
Let $(\tilde{\alpha}_{\dxp},\tilde{\beta}_{\dxp})$ be the solution of the same system 
\begin{align}
\tilde\alpha_{{\dxp}-1} &=\tilde\alpha_{\dxp} -\frac{i \idf^2}{2 \sin k} \Mpert_{\dxp} \big(\tilde\alpha_{\dxp} +\tilde\beta_{\dxp} e^{-2ik{\dxp}}\big),\\
\tilde\beta_{{\dxp}-1} &= \tilde\beta_{\dxp} +\frac{i \idf^2}{2 \sin k} \Mpert_{\dxp} \big(\tilde\alpha_{\dxp} e^{2ik{\dxp}}+\tilde\beta_{\dxp} \big),
\end{align}
for $0 < {\dxp} \leq {\Ld}$, $\dxp\in\Z,$ 
but with the  terminal conditions:
\begin{equation}
\tilde\alpha_{\Ld}=1, \quad \tilde\beta_{\Ld}=0 .
\end{equation}
Then, by linearity, we have
\begin{equation}
T = \frac{1}{\tilde{\alpha}_0},\quad R=\frac{\tilde{\beta}_0}{\tilde{\alpha}_0} .
\end{equation}
\begin{remark}
We can check that
\begin{equation}
|\tilde{\alpha}_{{\dxp}-1}|^2-|\tilde{\beta}_{{\dxp}-1}|^2 =
|\tilde{\alpha}_{{\dxp}}|^2-|\tilde{\beta}_{{\dxp}}|^2 ,
\end{equation}
which shows that $|\tilde{\alpha}_{{\dxp}}|^2-|\tilde{\beta}_{{\dxp}}|^2=1$ for all ${\dxp}\in\Z$, and therefore we get the energy conservation relation
\begin{equation}
|R|^2 +|T|^2=1.
\end{equation}
\end{remark}
\begin{remark}
If the variables $\Mpert_{\dxp}$ are independent and identically distributed with mean zero, then $(\tilde{\alpha}_{\dxp},\tilde{\beta}_{\dxp})$ is a martingale in the sense that: if we denote ${\mathcal F}_{\dxp}=\sigma( \Mpert_{{\dxp}'},  \, {\dxp}< {\dxp}'\leq {\Ld})$, then $(\tilde{\alpha}_{{\dxp}},\tilde{\beta}_{\dxp})$ is ${\mathcal F}_{\dxp}$-adapted and 
$\mathbb{E} \big[ (\tilde{\alpha}_{{\dxp}-1},\tilde{\beta}_{{\dxp}-1}) | {\mathcal F}_{\dxp} ]=  (\tilde{\alpha}_{{\dxp}},\tilde{\beta}_{\dxp})$.
\end{remark}

We apply the diffusion approximation theory (see Section \ref{app:adif}).
We find that
$\tau_{\Ld} := |T|^2$ behaves as a diffusion process as stated in Theorem~\ref{thm:squaretrans}.
This means that the probability density function of $\tau_{\Ld}$ satisfies the Fokker Planck equation
\begin{equation}\begin{split}
\label{eq:Fokker}
\partial_{\Ld} p_{\Ld}(\tau) = {\mathcal L}^* p_{\Ld} = \gamma \big[ \partial_\tau^2 \big( \tau^2(1-\tau) p_{\Ld}\big) +\partial_\tau\big( \tau^2 p_{\Ld}\big) \big],
\end{split}\end{equation}
starting from $p_{{\Ld}=0}(\tau)=\delta(\tau)$, where $\delta$ denotes the Dirac delta.
In particular we have for any $n\geq 1$:
\begin{equation}\begin{split}
\label{eq:momTa}
\mathbb{E} \big[ |T|^{2n} \big] &= \mathbb{E} \big[\tau_{\Ld}^n\big] = \int \tau^n p_{\Ld}(\tau) d\tau
\end{split}\end{equation}
which yields \eqref{eq:momT}.

\subsection{Scattering in non-matched medium}
\label{sec:proofnonmatched}
In this section we assume that 
\begin{align}
\Mpert_{\dxp} = \Delta_0 & \mbox{ for } {\dxp} \leq 0 ,\\
\Mpert_{\dxp} = \Delta_1 & \mbox{ for } {\dxp} > {\Ld} ,
\label{nonmatchmass}
\end{align}
$\dxp\in\Z$, and we introduce the wavenumbers ${\idk}_0$ and ${\idk}_1$ solutions of the dispersion relations
\begin{equation}
\label{eq:dispersionrelation:j2}
-(1+\Delta_j) \idf^2 =2 \cos {\idk}_j -2 - \Ks , \quad j=0,1.
\end{equation}
Here we assume the regime is propagative, i.e. the frequency $\idf$ is such that $(2+\Ks-\idf^2(1+\Delta_j))/(2) \in (-1,1)$ for $j=0,1$ so that there is a unique solution ${\idk}_j\in (0,\pi)$ to (\ref{eq:dispersionrelation:j}).

The wavenumber $k$ is still defined by (\ref{eq:dispersionrelation}).

We introduce 
\begin{align}
\alpha_{\dxp} &= \frac{e^{-ik{\dxp}}}{2i\sin k} \big( \hat{u}_{{\dxp}+1}-e^{-ik} \hat{u}_{\dxp} \big) ,\\ 
\beta_{\dxp} &= - \frac{e^{ik{\dxp}}}{2i\sin k} \big( \hat{u}_{{\dxp}+1}-e^{ik} \hat{u}_{\dxp} \big),
\end{align}
$\dxp\in\Z.$ We then have 
\begin{align}
\label{eq:bcnm1a}
\alpha_{\Ld} &= e^{i({\idk}_1-k){\Ld}} \frac{e^{i{\idk}_1}-e^{-ik}}{2i\sin k} T,\\
\beta_{\Ld} &= - e^{i({\idk}_1+k){\Ld}} \frac{e^{i{\idk}_1}-e^{ik}}{2i\sin k} T,
\label{eq:bcnm1b}
\end{align}
and
\begin{equation}
\hat{u}_{\dxp} = \alpha_{\dxp} e^{ik{\dxp}}+\beta_{\dxp} e^{-ik{\dxp}} ,
\end{equation}
for any ${\dxp}\in\Z$. The variables
$(\alpha_{\dxp},\beta_{\dxp})$ satisfy the system:
\begin{align}
\label{eq:sys2a}
\alpha_{{\dxp}-1} &=\alpha_{\dxp} -\frac{i \idf^2}{2 \sin k} \Mpert_{\dxp} \big(\alpha_{\dxp} +\beta_{\dxp} e^{-2ik{\dxp}}\big),\\
\beta_{{\dxp}-1} &=\beta_{\dxp} +\frac{i \idf^2}{2 \sin k} \Mpert_{\dxp} \big(\alpha_{\dxp} e^{2ik{\dxp}}+\beta_{\dxp} \big),
\label{eq:sys2b}
\end{align}
for ${\dxp} \leq {\Ld}$, with the terminal conditions (\ref{eq:bcnm1a}-\ref{eq:bcnm1b}) at ${\dxp}={\Ld}$. 

Moreover, we can express the reflection coefficient $R$ and the coefficient $\alpha_0$  in terms of $\beta_0  /\alpha_0$:
\begin{align}
\label{eq:expressRbeta0alpha0}
R &= - \frac{B+A \frac{ \beta_0}{\alpha_0} }
{ \overline{A} +\overline{B}  \frac{ \beta_0}{\alpha_0}},\\
\alpha_0 &= - \frac{2i \sin {\idk}_0}{  \overline{A} +\overline{B}  \frac{ \beta_0}{\alpha_0}},
\label{eq:expressalpha0beta0alpha0}
\end{align}
with
\begin{equation}
A=e^{i{\idk}_0}-e^{-ik} ,\qquad
B=e^{i{\idk}_0}-e^{ik}.
\end{equation}

\begin{proof}
We have $\hat{u}_0=1+R = \alpha_0+\beta_0$, therefore 
\begin{equation}
\alpha_0+\beta_0=1+R.
\end{equation}
We have $\hat{u}_{-1}=e^{-i{\idk}_0}+R e^{i{\idk}_0} = \alpha_{-1}e^{-ik}+\beta_{-1}e^{ik}$. We can express $(\alpha_{-1},\beta_{-1} )$ in terms of 
$(\alpha_0,\beta_0)$ by (\ref{eq:sys2a}-\ref{eq:sys2b}) evaluated at ${\dxp}=0$ (remember that $\Mpert_0=\Delta_0$), so that we get
\begin{equation}
\alpha_0 (e^{-ik}-{\idf^2}\Delta_0) +\beta_0(e^{ik}-{\idf^2}\Delta_0) =e^{-i{\idk}_0}+R e^{i{\idk}_0}  .
\end{equation}
By combining the two equations we can get the two relations
\begin{align}
\alpha_0 (e^{-ik}-e^{i{\idk}_0}-{\idf^2} \Delta_0) + \beta_0 (e^{ik}-e^{i{\idk}_0}- {\idf^2}\Delta_0)   \beta_0 = -2i \sin {\idk}_0,\\
\alpha_0 (- e^{-i{\idk}_0}+e^{-ik}-{\idf^2}\Delta_0) + \beta_0 (- e^{-i{\idk}_0}+e^{ik}-{\idf^2}\Delta_0) \beta_0 = 2i R \sin {\idk}_0,
\end{align}
which give:
\begin{align}
R &= \frac{(e^{-i{\idk}_0}-e^{-ik}+\frac{\idf^2}{K}\Delta_0) +(e^{-i{\idk}_0}-e^{ik}+{\idf^2}\Delta_0) \frac{ \beta_0}{\alpha_0} }
{ (e^{-ik}-e^{i{\idk}_0}-\frac{\idf^2}{K} \Delta_0) +(e^{ik}-e^{i{\idk}_0}-{\idf^2}\Delta_0) \frac{ \beta_0}{\alpha_0}},\\
\alpha_0 &=  \frac{-2i \sin {\idk}_0}{ (e^{-ik}-e^{i{\idk}_0}-{\idf^2} \Delta_0) +(e^{ik}-e^{i{\idk}_0}-{\idf^2}\Delta_0) \frac{ \beta_0}{\alpha_0}}.
\end{align}
We then get the desired result by remarking that $\frac{\Delta_0 \idf^2}{2K} = \cos k -\cos {\idk}_0$, which gives $ e^{-ik}-e^{i{\idk}_0}-{\idf^2} \Delta_0 = -e^{ik}+e^{-i{\idk}_0} = \overline{A}$ and $e^{ik}-e^{i{\idk}_0}-{\idf^2}\Delta_0=e^{-i{\idk}_0}-e^{-ik} =\overline{B}$. 
\end{proof}

Let $(\tilde{\alpha}_{\dxp},\tilde{\beta}_{\dxp})$ be the solution of the same system (\ref{eq:sys2a}-\ref{eq:sys2b})
\begin{align}
\tilde\alpha_{{\dxp}-1} &=\tilde\alpha_{\dxp} -\frac{i \idf^2}{2 \sin k} \Mpert_{\dxp} \big(\tilde\alpha_{\dxp} +\tilde\beta_{\dxp} e^{-2ik{\dxp}}\big),\\
\tilde\beta_{{\dxp}-1} &= \tilde\beta_{\dxp} +\frac{i \idf^2}{2 \sin k} \Mpert_{\dxp} \big(\tilde\alpha_{\dxp} e^{2ik{\dxp}}+\tilde\beta_{\dxp} \big),
\end{align}
for $0 < {\dxp} \leq {\Ld}, \dxp\in\Z$, 
but with the  terminal conditions at ${\dxp}={\Ld}$:
\begin{equation}
\tilde\alpha_{\Ld}=e^{i({\idk}_1-k){\Ld}} \frac{e^{i{\idk}_1}-e^{-ik}}{2i\sin k} , \quad 
\tilde\beta_{\Ld}= - e^{i({\idk}_1+k){\Ld}} \frac{e^{i{\idk}_1}-e^{ik}}{2i\sin k}.
\end{equation}
Then, by linearity, we have $\tilde\beta_0/\tilde\alpha_0=\beta_0/\alpha_0$ so we get from (\ref{eq:expressRbeta0alpha0}):
\begin{align}
R= -\frac{A}{\overline{A}} \frac{\frac{B}{A}+\tilde{R}}{1+\frac{\overline{B}}{\overline{A}} \tilde{R}},
\end{align}
with
\begin{align}
\tilde{R}=\frac{\tilde{\beta}_0}{\tilde{\alpha}_0} .
\end{align}
By linearity, we have $\alpha_0=\tilde{\alpha}_0 T$ so that we get from (\ref{eq:expressalpha0beta0alpha0}):
\begin{equation}
T = - \frac{2i \sin {\idk}_0}{\overline{A}}\frac{1}{1+\frac{\overline{B}}{\overline{A}} \tilde{R}} \tilde{T},
\end{equation}
with
\begin{equation}
\tilde{T} = \frac{1}{\tilde{\alpha}_0}.
\end{equation}

We can check that
\begin{equation}
|\tilde{\alpha}_{{\dxp}-1}|^2-|\tilde{\beta}_{{\dxp}-1}|^2 =
|\tilde{\alpha}_{{\dxp}}|^2-|\tilde{\beta}_{{\dxp}}|^2 ,
\end{equation}
which shows that $|\tilde{\alpha}_{{\dxp}}|^2-|\tilde{\beta}_{{\dxp}}|^2 = |\tilde{\alpha}_{{\Ld}}|^2-|\tilde{\beta}_{{\Ld}}|^2  =\frac{\sin {\idk}_1}{\sin k}$ for all ${\dxp}$, and therefore we get the energy conservation relation
\begin{equation}
|R|^2 + \frac{\sin {\idk}_1}{\sin k} |T|^2=1.
\end{equation}

In case ${\Mpt}_0\ne{\Mpt}_1, {\Mpt}_0{\Mpt}_1=0$, we have
\begin{align}
|T|^2 &=\frac{\sin k}{\sin {\idk}_1} |\tilde{T}|^2 ,
\end{align}
where $ |\tilde{T}|^{2}$ behaves as the diffusion process $\tau_{\Ld}$ with the infinitesimal generator
(\ref{eq:defL}) starting from 
\begin{equation}
\tau_0 = \frac{2\sin k \sin {\idk}_1}{1-\cos(k+{\idk}_1)}.
\end{equation}
We get the following representation of the probability density function of $|\tilde{T}|^{2}$:
\begin{equation}\begin{split}
p_{\Ld}(\tau) =\frac{2}{\tau^2} e^{  - \frac{\gamma {\Ld}}{4} }
\int_0^\infty s \tanh(\pi s) P_{-1/2+is} \Big(\frac{2}{\tau}-1\Big)\\
\times P_{-1/2+is} \Big(\frac{2}{\tau_0}-1\Big) e^{  -s^2 \gamma {\Ld} } ds,
\end{split}\end{equation}
where $P_{-1/2+is}(\eta)$, $\eta\geq 1$, $s \geq 0$
is the Legendre function of the first
kind, which is the solution of
\begin{equation}\begin{split}
\label{eq:legendre}
\frac{d}{d\eta}(\eta^2-1)\frac{d}{d\eta}P_{-1/2+is}(\eta) =
- \left( s^2+ \frac{1}{4} \right) P_{-1/2+is}(\eta)\, ,
\end{split}\end{equation}
starting from $P_{-1/2+i s}(1)=1$.
It has the integral representation \eqref{Phalfeta}.
In particular, we have \eqref{nonmatchedcase2} and \eqref{nonmatchedcase2nn}.

\begin{remark}
It is possible to address the case where ${\idf}$ is outside the common pass band on the right half-space, i.e. $(2+\Ks-{\idf}^2(1+{\Mpt}_1))/(2) \not\in (-1,1)$. 
If $-1<{\Mpt}_1<0$ is such that $(2+\Ks-{\idf}^2(1+{\Mpt}_j))/(2)>1$, then the wave has the form
$
\hat{u}_{\dxp} = T e^{-{\idk}_1 {\dxp}}
$
for ${\dxp}>{\Ld}$ instead of (\ref{eq:formuyafetrL}), where ${\idk}_1$ is given by
\begin{equation}\begin{split}
\cosh({\idk}_1) = \frac{2+\Ks-{\idf}^2(1+{\Mpt}_1)}{2}.
\end{split}\end{equation}
We can then proceed as above and find that, instead of (\ref{eq:bcnm1a}-\ref{eq:bcnm1b}), we have
\begin{equation}\begin{split}
\alpha_{\Ld} &= e^{(-{\idk}_1-ik){\Ld}} \frac{e^{-{\idk}_1}-e^{-ik}}{2i\sin k} T,\\
\beta_{\Ld} &= - e^{(-{\idk}_1+ik){\Ld}} \frac{e^{-{\idk}_1}-e^{ik}}{2i\sin k} T.
\end{split}\end{equation}
This implies that $|\tilde{\alpha}_{\dxp}|^2-|\tilde{\beta}_{\dxp}|^2=|\tilde{\alpha}_{\Ld}|-|\tilde{\beta}_{\Ld}|^2=0$ for all $0\leq {\dxp} \leq {\Ld}, \dxp\in\Z$, and therefore $|\tilde{R}|=1$ and
$|R|=1$. The wave is totally reflected, the random section only changes the phase of the reflected component compared to the case without random perturbation:
\begin{equation}\begin{split}
\hat{u}_{\dxp} = e^{i{\idk}_0 {\dxp}}+R e^{-i{\idk}_0{\dxp}} \mbox{ for } {\dxp} \leq 0, \dxp\in\Z.
\end{split}\end{equation}

If ${\Mpt}_1>0$ is such that $Q:=(2+\Ks-{\idf}^2(1+{\Mpt}_j))/(2) <-1$, then the wave has the form
$
\hat{u}_{\dxp} = T \rho^{\dxp}
$
for ${\dxp}>{\Ld}$ instead of (\ref{eq:formuyafetrL}), where $\rho$ is given by
$
\rho=Q+\sqrt{Q^2-1}
$
(it is the unique solution in $(-1,1)$ of the dispersion relation $-2Q+\rho+\rho^{-1}=0$). We obtain the same conclusion: the reflection coefficient $R$ has modulus one.
\end{remark}

\section{Proofs of Theorems 
\ref{thmTDmeanks0},
\ref{thmTDfrontks0},
\ref{thmTDmeanks},
\ref{thmTDfrontks}
}
\label{sec:prooftd}
Here we consider that $\dxp_0=0$ and that in the section $[1,{\Ld}] \cap \mathbb{Z}$ the variables $\Mpert_{\dxp}$ are independent and identically distributed with mean zero and variance $\sigma^2$.
The transmitted wave for $\dxp>{\Ld}$ has the form
\begin{equation}
{\su}_{{\dxp}}(t) = \frac{1}{2\pi} \int_{-\infty}^{+\infty}\hat{\su}_{{\dxp}}(\fq) e^{-i \fq t} d\fq , \qquad t \geq 0, \quad \dxp \in \mathbb{Z},
\end{equation}
where the Fourier components are given by 
\begin{equation}
\hat{\su}_{{\dxp}}(\fq) =  
\hat{c}(\fq) T(\fq) \cos \big( k(\fq) (\dxp-\dxp_0) \big), 
\end{equation}
with $k(\fq)$ and $\hat{c}(\fq)$ given by (\ref{eq:defkc}) when $\Ks=0$ and by (\ref{eq:defkomegas}) for $\Ks \geq 0$.
The statistics of the transmission coefficient $T(\fq)$ at a fixed frequency has been studied in the previous section. We have 
\begin{equation}
\mathbb{E} \big[ T(\fq) \big]=e^{ -\gamma(\fq) {\Ld} },
\end{equation}
\begin{equation}
\mathbb{E}\big[|T(\fq) |^2\big]=e^{ -\frac{1}{4}\gamma(\fq) {\Ld} }
\int_0^\infty e^{- \gamma(\fq) {\Ld} s^2} \frac{2 \pi s \sinh(\pi s)}{\cosh^2(\pi s)}  ds,
\end{equation}
where $\gamma(\fq)$ is given by (\ref{eq:defgamma2}).
In order to characterize the time-dependent field, we need to characterize the statistics of the vector $(T(\fq_j))_{j=1}^n$ for any set of distinct frequencies $(\fq_j)_{j=1}^n$, more exactly, we need to characterize the moments
\begin{equation}
\mathbb{E}\Big[ \prod_{j=1}^n T(\fq_j) \Big] 
\end{equation}
because they in turn characterize all moments of the field
\begin{equation}\begin{split}
\mathbb{E}\big[ {\su}_{{\dxp}}(t)^n \big] = \frac{1}{(2\pi)^n} \int_{-\infty}^{+\infty} \cdots\int_{-\infty}^{+\infty} 
\mathbb{E} \Big[  \prod_{j=1}^n T(\fq_j) \Big] 
\prod_{j=1}^n \hat{c}(\fq_j)  \cos \big( k(\fq_j) (\dxp-\dxp_0) \big)\\
e^{-i \sum_{j=1}^n \fq_j t} 
d\fq_1\cdots d\fq_n.
\end{split}\end{equation}
Proceeding as in \cite[Chapter 8]{book}, we can show that,  for any set of distinct frequencies $(\fq_j)_{j=1}^n$, 
$(T(\fq_j))_{j=1}^n$ has the martingale representation 
\begin{equation}
T(\fq_j) = M(\fq_j) \tilde{T}(\fq_j),
\end{equation} 
where 
\begin{equation}
\tilde{T}(\fq) = \exp \big( i \sqrt{\gamma(\fq)} W_{\Ld}- \gamma(\fq){\Ld}/2 \big), 
\end{equation}
$W_{\Ld} \sim {\mathcal N}(0,{\Ld})$, 
and $M(\fq_j)$ are independent complex martingales (and independent of $W_{\Ld}$) with mean one.
As a consequence,
\begin{equation}
\mathbb{E} \Big[  \prod_{j=1}^n T(\fq_j) \Big] 
=
\mathbb{E}\Big[ \prod_{j=1}^n \tilde{T}(\fq_j) \Big] ,
\end{equation}
and ${\su}_{{\dxp}}(t)$ can be written as (more exactly, it has the same moments as)
\begin{equation}
{\su}_{{\dxp}}(t) = \frac{1}{2\pi} \int_{-\infty}^{+\infty} \tilde{T}(\fq)  
\hat{c}(\fq) \cos \big( k(\fq) (\dxp-\dxp_0) \big)  e^{-i \fq t} d\fq .
\end{equation}

\section*{Acknowledgments}
This work was started during the stay (in March 2023) of both authors at the Isaac Newton Institute (INI) for Mathematical Sciences, Cambridge. The authors would like to thank INI for support and hospitality during the programme -- `Mathematical theory and applications of multiple wave scattering' (MWS) where work on this paper was undertaken. This work was supported by EPSRC grant no EP/R014604/1. A part of the work of BLS, for the same visit to INI, was partially supported by a grant from the Simons Foundation.
BLS thanks P. Martin for stimulating discussions during MWS regarding solution in \cite{Schrodinger}.

\appendix

\section{Exact solutions}

\subsection{Matched medium}
\label{Tfullmatchedapp}
Using the Green's function for one-dimensional model,
the scattered field is given by
\begin{equation}\begin{split}
\hat{\su}_{\dxp}
&=-\idf^2{\sC}({\dxk})\sum_{q=1}^{{{\Nd}}}{{\Mpert}}_q(\hat{\su}_{q}+\zxk^q)e^{{i{\dxk}}{|\dxp-q|}}, \quad \dxp\in\Z,
\label{matchedux}
\end{split}\end{equation}
with
\begin{equation}\begin{split}
{\sC}({\dxk}):=\frac{1}{2i\sin\dxk},\quad \zxk:=e^{{i{\dxk}}}\in\mathbb{C}.
\label{defC2}
\end{split}\end{equation}
The reduced set of equations is therefore
\begin{equation}\begin{split}
\hat{\su}_{p}
&=-\idf^2{\sC}({\dxk})\sum_{q=1}^{{{\Nd}}}{{\Mpert}}_q{\zxk}^{|p-q|}(\hat{\su}_{q}+\zxk^q),\quad p=[1,{{\Nd}}]\cap\Z.
\label{matchedup}
\end{split}\end{equation}
Using the definitions
\begin{align}
\mbf{u}:=(\hat{\su}_{1}, \hat{\su}_{2}, \dotsc, \hat{\su}_{\Nd})^T\in\mathbb{C}^{\Nd},\label{uvec}
\end{align}
\begin{subequations}
\begin{align}
\mbf{z}_{{\Nd}}({\mathfrak{z}}):=(1, {\mathfrak{z}}, {\mathfrak{z}}^2, \dotsc, {\mathfrak{z}}^{{{\Nd}}-1})^T\in\mathbb{C}^{\Nd}, \quad {\mathfrak{z}}\in\mathbb{C},
\label{zNdef}\\
\mbf{D}({\Mpert}):=-\idf^2\text{diag}({{\Mpert}}_1,{{\Mpert}}_2,\dotsc,{{\Mpert}}_{{\Nd}})\in\mathbb{C}^{\Nd\times\Nd},\label{Dmat}\\
\mbf{I}:=\text{diag}(1, \dotsc, 1)\in\mathbb{C}^{\Nd\times\Nd}\quad\text{ (identity matrix)},\label{Imat}
\end{align}
\label{matcheduTD}
\end{subequations}
\begin{equation}
\mbf{T}({\zxk}):={\sC}({\dxk})\text{Toeplitz}(\mbf{z}_{{\Nd}}({\zxk})^T)\in\mathbb{C}^{\Nd\times\Nd},\label{Tmat}
\end{equation}
and \eqref{zNdef}, the question of scattered field can therefore be answered formally by \eqref{matchedup} as
\begin{equation}
\mbf{u}=(\mbf{I}-\mbf{T}({\zxk})\mbf{D}({\Mpert}))^{-1} \mbf{T}({\zxk})\mbf{D}({\Mpert}){\zxk}\mbf{z}_{{\Nd}}({\zxk}),
\label{matchedusol}
\end{equation}
where %$\mbf{I}$ is the ${\Nd\times\Nd}$ identity matrix, 
$\zxk$ is defined in \eqref{defC2}${}_2$ and ${\sC}$ in \eqref{defC2}${}_1$.
For $\dxp>{\Nd}, \dxp\in\Z,$ according to \eqref{matchedux},
\begin{equation}\begin{split}
\hat{\su}_{\dxp}&=-\idf^2{\sC}({\dxk})\sum_{q=1}^{{{\Nd}}}{{\Mpert}}_q{\zxk}^{\dxp-q}(\hat{\su}_{q}+\zxk^q),
\end{split}\end{equation}
so that, using the definitions \eqref{defC2},\eqref{zNdef}, \eqref{matcheduTD} and the expression \eqref{matchedusol}, the transmission coefficient $T$ in
\begin{equation}\begin{split}
\hat{\su}_{\dxp,0}+\zxk^{\dxp}
&=T\zxk^{\dxp},
\label{matchedTasymp}
\end{split}\end{equation}
is
\begin{equation}\begin{split}
T=1+{\sC}({\dxk})\mbf{z}_{{\Nd}}({\zxk}^{-1})\cdot\mbf{D}({\Mpert})(\mbf{I}-\mbf{T}({\zxk})\mbf{D}({\Mpert}))^{-1}\mbf{z}_{{\Nd}}({\zxk})\in\mathbb{C}.
\label{Tfullmatched}
\end{split}\end{equation}
In the absence of mass perturbation on $[1,\Ld]\cap\Z$, it is clear that $\mbf{D}({\Mpert})=\mbf{0}$ implies $T=1$, as expected.

\subsection{Non-matched medium}
%Recall the definitions \eqref{zNdef} and \eqref{defC2}${}_2$ in the following.
Recall
\eqref{nonmatchmass}
and
\eqref{eq:dispersionrelation:j2}
regarding the definitions of ${\Mpt}_0, {\Mpt}_1$ and $\idk_0, \idk_1$, respectively.

\subsubsection{Case 1: {${\Mpt}_0\ne0, {\Mpt}_1=0$}{}}
Let 
\begin{equation}
\mbf{u}^\ast:=(1+{\sC}_{2}({\idk}_{0},k))(\mbf{I}-\mbf{T}\mbf{D}({\Mpert}))^{-1}\zxk\mbf{z}_{{\Nd}}(\zxk)\in\mathbb{C}^{\Nd},
\label{nonmatchedustar}
\end{equation}
using
the definitions \eqref{matcheduTD} and \eqref{defC2}${}_2$,
%$\zxk$ is defined in \eqref{defC2}${}_2$,
and
\begin{equation}
\mbf{T}=[T_{pq}]_{p, q=1}^{{\Nd}}:=[G_{p,q}]_{p, q=1}^{{\Nd}}\in\mathbb{C}^{\Nd\times\Nd},
\end{equation}
with
\begin{equation}\begin{split}
{G}_{{p}, {q}}:=
\begin{cases}
e^{-i {\idk}_0({p}-1)}e^{i {\idk}_0}{{\sC}}_{3}({\idk}_{0},k)e^{ik({q}-1)}, {p}<1\\
\left(e^{ik{({p}-1)}}{{\sC}}_1({\idk}_0,k)e^{2i {k}}+e^{-ik{({p}-1)}}{{\sC}}(k)\right)e^{ik{({q}-1)}}, 1\le {p}\le {q}\\
e^{ik{({p}-1)}} \left(e^{ik{({q}-1)}}{{\sC}}_1({\idk}_0,k)e^{2i {k}}+e^{-ik{({q}-1)}}{{\sC}}(k)\right), {p}> q,
\end{cases}
\label{Gij}
\end{split}\end{equation}
$p,q\in\Z,$ where, in addition to \eqref{defC2}${}_1$, we employ the definitions
\begin{equation}\begin{split}
{{\sC}}_1({\idk}_0,k):={{\sC}}(k){\sC}_{2}(k,{\idk}_{0}),  
{\sC}_{2}({\idk}_{0},k):=-\frac{1-e^{i (k-{\idk}_{0})}}{1-e^{i (k+{\idk}_{0})}},
{{\sC}}_{3}({\idk}_{0},k):=-\frac{e^{i (k-{\idk}_{0})}}{1-e^{i({\idk}_{0}+k)}}.
\label{C1C2C3}
\end{split}\end{equation}

For $\dxp>{\Nd}, \dxp\in\Z,$ as $\dxp\to+\infty$,
using \eqref{Gij}${}_3$,
the transmission coefficient in \eqref{matchedTasymp} is found to be
\begin{equation}\begin{split}
T=1+{\sC}_{2}({\idk}_{0},k)+T_1+T_2\in\mathbb{C},
\label{nonmatchexactT}
\end{split}\end{equation}
with
\begin{equation}\begin{split}
{T}_1:={{\sC}}_1({\idk}_0,k)
\zxk\mbf{z}_{{\Nd}}(\zxk)\cdot\mbf{D}({\Mpert})\mbf{u}^\ast,\quad
{T}_2:={{\sC}}(k)
\zxk^{-1}\mbf{z}_{{\Nd}}(\zxk^{-1})\cdot\mbf{D}({\Mpert})\mbf{u}^\ast,
\label{nonmatchexactT1both}
\end{split}\end{equation}
where we also use the definitions \eqref{defC2}, \eqref{matcheduTD}, \eqref{nonmatchedustar}, \eqref{C1C2C3}.

\begin{remark}
As ${\Mpt}_0\to0$, we get
${{\sC}}_1({\idk}_0,k)={{\sC}}(k){\sC}_{2}(k,{\idk}_{0})\to0,$
${{\sC}}_{3}({\idk}_{0},k)\to-{1}/{ \left(1-e^{2ik}\right)}, {\sC}_{2}({\idk}_{0},k)\to0.$
Hence,
${G}_{{p}, {q}}\to{{\sC}}(k)e^{ik|p-{q}|},$
and \eqref{Tfullmatched} follows from \eqref{nonmatchexactT} after simplifying \eqref{nonmatchedustar}.
\end{remark}

\subsubsection{Case 2: {${\Mpt}_0=0, {\Mpt}_1\ne0$}{}}
%Recall \eqref{zNdef}.
Let 
\begin{equation}
\mbf{u}^\ast:=(\mbf{I}-\mbf{T}\mbf{D}({\Mpert}))^{-1}({\zxk}\mbf{z}_{{\Nd}}(\zxk)+{\sC}_{2}({k},{\idk}_1){\zxk}^{2\Ld-1}\mbf{z}_{{\Nd}}(\zxk^{-1}))\in\mathbb{C}^{\Nd},
\label{nonmatched2ustar}
\end{equation}
using
the definitions \eqref{matcheduTD} and \eqref{defC2}, \eqref{C1C2C3},
and
\begin{equation}
\mbf{T}=[T_{pq}]_{p, q=1}^{\Nd}:=[G_{p,q}]_{p, q=1}^{\Nd}\in\mathbb{C}^{\Nd\times\Nd},
\end{equation}
with
\begin{equation}\begin{split}
{G}_{{p}, {q}}:=
\begin{cases}
e^{-i {k} {(q-\Ld-1)}}e^{i{k}}{{\sC}}_{3}({k},{\idk}_{1})e^{i{\idk}_{1} {(p-\Ld-1)}}, {p}>\Ld,\\
e^{-i {k} {(q-\Ld-1)}}\left(e^{-i{k} {(p-\Ld-1)}}{{\sC}}_1({\idk}_1,{k})+e^{i{k} {(p-\Ld-1)}}{{\sC}}({k})\right), {q}\le {p}\le\Ld,\\
e^{-i {k} {(p-\Ld-1)}}\left(e^{-i{k} {(q-\Ld-1)}}{{\sC}}_1({\idk}_1,{k})+e^{i{k} {(q-\Ld-1)}}{{\sC}}({k})\right), {p}<{q},
\end{cases}
\label{Gij2}
\end{split}\end{equation}
$p,q\in\Z$.

For $\dxp>\Ld, \dxp\in\Z,$ as $\dxp\to+\infty$,
using \eqref{Gij2}${}_1$,
the transmission coefficient in
\begin{equation}\begin{split}
\hat{\su}_{\dxp,0}+\zxk^{\dxp}&=T e^{i{\idk}_{1}\dxp},
\end{split}\end{equation}
is obtained as
\begin{equation}\begin{split}
T={\zxk}^{\Ld}e^{-i\idk_{1}\Ld}(1+{\sC}_{2}({k},{\idk}_{1})+T_3)\in\mathbb{C},
\label{nonmatchexactT2a}
\end{split}\end{equation}
with
\begin{equation}\begin{split}
{T}_3&:=e^{-i{\idk}_{1}}{\sC}_{3}({k},{\idk}_{1})\zxk\mbf{z}_{{\Nd}}(\zxk^{-1})\cdot\mbf{D}({\Mpert})\mbf{u}^\ast,
\label{nonmatchexactT12}
\end{split}\end{equation}
where the definitions \eqref{defC2}, \eqref{matcheduTD}, \eqref{nonmatched2ustar}, \eqref{C1C2C3} are employed.

\begin{remark}
As ${\Mpt}_1\to0$, we get
${{\sC}}_1({\idk}_1,{k})={{\sC}}({k}){\sC}_{2}({k},{\idk}_{1})\to0,$
${{\sC}}_{3}({k},{\idk}_1)\to-{1}/{ \left(1-e^{2i{k}}\right)}, {\sC}_{2}({k},{\idk}_1)\to0.$
Hence,
${G}_{{p}, {q}}
\to
{{\sC}}(k)e^{ik|p-{q}|},$
and \eqref{Tfullmatched} follows from \eqref{nonmatchexactT2a}
after simplifying \eqref{nonmatched2ustar}.
\end{remark}

\begin{remark}
Due to the definitions \eqref{nonmatchexactT1both} and \eqref{nonmatchexactT12},
in the absence of mass perturbation on $[1,\Ld]\cap\Z$, it is clear that $\mbf{D}({\Mpert})=\mbf{0}$ so that \eqref{nonmatchexactT} and \eqref{nonmatchexactT2a} lead to $|T|=|1+{\sC}_{2}({\idk_{0}},{\idk}_{1})|$. Upon simplifying $|T|^2$ further it is easy find that
\begin{equation}\begin{split}
|T|^2&
=\frac{1-\cos(2\idk_{0})}{1-\cos(\idk_{0}+\idk_{1})}.
\label{nonmatchexactTperf}
\end{split}\end{equation}
Indeed, $|T|=1$ in matched case without mass perturbation as $\idk_{0}=\idk_{1}$.
\label{remarkTmisperf}
\end{remark}

\section{Diffusion-approximation}
\label{app:adif}
Here we give a few elements to the proof of the convergence of the process $|T|^2$ to the diffusion Markov process with the generator given by Eq.~(\ref{eq:defL}).
The proof consists in showing that $(\tilde{\alpha}_{\dxp},\tilde{\beta}_{\dxp})$ converges to a diffusion Markov process,
that $\tilde{\beta}_{\dxp}/\tilde{\alpha}_{\dxp}$ is itself a diffusion Markov process, that $|\tilde{\beta}_{\dxp}/\tilde{\alpha}_{\dxp}|^2$ is itself a diffusion Markov process, and then compute the moments of $|T|^2=1-|\tilde{\beta}_0/\tilde{\alpha}_0|^2$.

The first diffusion-approximation theorem is the following one:\\
\begin{prop}
Let $X_j$ be a $\mathbb{R}^d$-valued random sequence 
solution of 
\begin{equation}\begin{split}
X_{j+1} = X_j+ \epsilon Y_j F(X_j) ,
\end{split}\end{equation}
where $Y_j$ are independent and identically distributed with mean zero and variance $\sigma^2$ and $F:\mathbb{R}^d\to\mathbb{R}^d$ is smooth. 
Let $X^\epsilon(z) = X_{[z\epsilon^{-2}]}$, where $[\cdot]$ stands for the integer part.
As $\epsilon \to 0$ the process $X^\epsilon$ converges (in the space of the cadlag functions) to the diffusion process $\bar{X}$ with the infinitesimal generator
\begin{equation}\begin{split}
{\mathcal L} = \frac{\sigma^2}{2} \sum_{i,i'=1}^d F_i(x) F_{i'}(x) \partial^2_{x_ix_{i'}} ,
\end{split}\end{equation}
or, equivalently, solution of the stochastic differential equation
\begin{equation}\begin{split}
d\bar{X} = \sigma F(\bar{X}) dW(z) ,
\end{split}\end{equation}
where $W$ is a Brownian motion and the stochastic integral is It\^o.
\end{prop}

We can use the diffusion-approximation theorem 6.1 in \cite[Chapter 6]{book} to prove this result. We need to pay attention to the fact that we here deal with a discrete system, not a continuous one. We need to introduce the correct approximation of $X^\epsilon$ which is:
\begin{equation}\begin{split}
\frac{d \tilde{X}^\epsilon(z)}{dz} = \frac{1}{\eps} Y_{[\epsilon^{-2} z]} F ( \tilde{X}^\epsilon ) -\frac{1}{2} \sum_{i=1}^d \partial_{x_i} F( \tilde{X}^\epsilon) F_i ( \tilde{X}^\epsilon) Y_{[\epsilon^{-2} z]}^2.
\end{split}\end{equation}
Then we can show that $ \tilde{X}^\epsilon(z) - X^\epsilon(z)$ converges to zero and by \cite[Theorem 6.1]{book} that $\tilde{X}^\epsilon(z)$ converges to a diffusion process with the infinitesimal generator
\begin{equation}\begin{split}
{\mathcal L} = \frac{\sigma^2}{2} \sum_{i,i'=1}^d F_i(x) \partial_{x_i}\big( F_{i'}(x) \partial_{x_{i'}} \big) 
-\frac{\sigma^2 }{2} \sum_{i,i'=1}^d \partial_{x_i} F_{i'} ( x) F_i (x) \partial_{x_{i'}} \\
= \frac{\sigma^2}{2} \sum_{i,i'=1}^d F_i(x) F_{i'}(x) \partial^2_{x_ix_{i'}}.
\end{split}\end{equation}
In fact, the result can also be obtained from \cite[Chapter 7, Corollary 4.2]{ethier} which addresses directly discrete systems. In the framework of \cite{ethier}, $\epsilon=1/\sqrt{n}$ and the transition function is 
\begin{equation}\begin{split}
\mu_n (x,\Gamma) = \mathbb{P}\big( X_{j+1} \in \Gamma | X_j=x\big) = 
\mathbb{P} \big( x+ \sqrt{n}^{-1} Y F(x) \in \Gamma\big) ,
\end{split}\end{equation}
which is such that 
\begin{align}
n \int (\tilde{x}-x) \mu_n(x,d\tilde{x}) &= 0 , \\
n \int (\tilde{x}-x)(\tilde{x}-x)^T \mu_n(x,d\tilde{x}) &= \sigma^2 F(x) F(x)^t.
\end{align}

We need in fact a version of the diffusion-approximation theorem with a periodic component.
The diffusion-approximation theorem that we need is the following one:
\begin{prop}
Let $X_j$ be a $\mathbb{R}^d$-valued random sequence 
solution of 
\begin{equation}\begin{split}
X_{j+1} = X_j+ \epsilon Y_j F(X_j,j) ,
\end{split}\end{equation}
where $Y_j$ are independent and identically distributed with mean zero and variance $\sigma^2$ and $F$ is smooth with respect to its first entry and periodic with respect to its second entry. 
Let $X^\epsilon(z) = X_{[z\epsilon^{-2}]}$, where $[\cdot]$ stands for the integer part.
As $\epsilon \to 0$ the process $X^\epsilon$ converges to the diffusion process $\bar{X}$ with the infinitesimal generator
\begin{equation}\begin{split}
{\mathcal L} = \frac{\sigma^2}{2} \sum_{i,i'=1}^d \left< F_i(x,\cdot) F_{i'}(x,\cdot) \right> \partial^2_{x_ix_{i'}} ,
\end{split}\end{equation}
where $\left< \cdot \right>$ is an average over the periodic component. Equivalently, $\bar{X}$ is solution of the stochastic differential equation
\begin{equation}\begin{split}
d\bar{X} = \sum_{j=1}^{d'} \sigma \bar{F}^{(j)} (\bar{X}) dW^{(j)}(z) ,
\end{split}\end{equation}
where the $W^{(j)}$ are independent Brownian motions, the stochastic integrals are It\^o and we have identified $d'\geq 1$ functions $F^{(j)}:\mathbb{R}^d\to\mathbb{R}^d$ such that 
\begin{equation}
\sum_{j=1}^{d'} \bar{F}^{(j)}(x) \bar{F}^{(j)}(x)^T = \left< F(x,\cdot) F(x,\cdot)^T \right>.
\end{equation}
\end{prop}

We can use the diffusion-approximation theorem 6.5 in \cite[Chapter 6]{book} to prove this result.
When $L=[z_0/\sigma^2]$ for some $z_0>0$,
the application of this theorem (proceeding as in \cite[Chapter 7]{book}) gives that $X^\sigma(z) = (\tilde{\alpha}_{[z_0/\sigma^2-z/\sigma^2]},
\tilde{\beta}_{[z_0/\sigma^2-z/\sigma^2]})^T$ converges as $\sigma \to 0$ to a diffusion Markov process $\bar{X}(z)$, $z\in [0,z_0]$, solution of the stochastic differential equation
\begin{equation}\begin{split}
d\bar{X} = 
\frac{i {\fq}^2 \sigma}{2  \sin k } 
\begin{pmatrix} 1 &0\\ 0 & -1\end{pmatrix} \bar{X} dW_0(z)
- \frac{ {\fq}^2 \sigma}{2 \sqrt{2}  \sin k } 
\begin{pmatrix} 0 &1\\ 1 & 0\end{pmatrix} \bar{X} dW_1(z)\\
- \frac{i {\fq}^2 \sigma}{2 \sqrt{2}  \sin k } 
\begin{pmatrix}0 &1\\ -1 & 0\end{pmatrix} \bar{X} d\tilde{W}_1(z) ,
\end{split}\end{equation}
where $W_0$, $W_1$ and $\tilde{W}_1$ are independent Brownian motions.

Denoting $\tilde{R}^\sigma(z) = \tilde{\beta}_{[z_0/\sigma^2-z/\sigma^2]}/
\tilde{\alpha}_{[z_0/\sigma^2-z/\sigma^2]}$, we deduce that $\tilde{R}^\sigma(z)$ converges as $\sigma \to 0$ to a diffusion Markov process $\tilde{R}(z)$, $z\in [0,z_0]$, solution of the stochastic differential equation
\begin{equation}\begin{split}
d\tilde{R} = 
- \frac{i {\fq}^2 \sigma}{  \sin k } \tilde{R} dW_0(z)
+ \frac{ {\fq}^2 \sigma}{2 \sqrt{2}  \sin k } 
\big(\tilde{R}^2-1\big) dW_1(z)\\
+ \frac{i {\fq}^2 \sigma}{2 \sqrt{2}  \sin k } 
\big(\tilde{R}^2+1\big) d\tilde{W}_1(z)
-\frac{ 3{\fq}^4 \sigma^2}{4 \sin^2 k } \tilde{R} dz.
\end{split}\end{equation}


\begin{thebibliography}{99}

\bibitem{Amir}
Amir, A., Oreg, Y., Imry, Y. (2018). Thermal conductivity in 1d: Disorder-induced transition from anomalous to normal scaling. Europhysics letters, 124(1), 16001

\bibitem{Anderson}
Anderson, P. W. (1958). Absence of diffusion in certain random lattices. Physical review, 109(5), 1492

\bibitem{balgu15}
Bal, G., Gu, Y. (2015).
Limiting models for equations with large random potential; a review. 
Communication in Mathematical Sciences 13, 729--748. 

\bibitem{Basile2}
Basile, G., Bernardin, C., Jara, M., Komorowski, T., Olla, S. (2016). Thermal conductivity in harmonic lattices with random collisions. Thermal transport in low dimensions: from statistical physics to nanoscale heat transfer, 215--237

\bibitem{Basile2014}
Basile, G., Olla, S. (2014). Energy diffusion in harmonic system with conservative noise. Journal of Statistical Physics, 155(6), 1126--1142

\bibitem{Basile}
Basile, G., Olla, S., Spohn, H. (2010). Energy transport in stochastically perturbed lattice dynamics. Archive for rational mechanics and analysis, 195(1), 171--203

\bibitem{Bloch}
Bloch, F. (1929). \"{U}ber die Quantenmechanik der Elektronen in Kristallgittern. Z. Physik 52, 555--600

\bibitem{BornHuang1985}
Born, M., Huang, K. (1985). Dynamical theory of crystal lattices. Oxford Classic Texts in the Physical Sciences. The Clarendon Press, Oxford University Press, New York

\bibitem{BornHuang1912}
Born, M., von Karman, T. (1912). On fluctuations in spatial grids. Phys. Z. 13, 297--309

\bibitem{Brillouin}
Brillouin, L. (1953).
{\em Wave propagation in periodic structures; electric filters and crystal lattices}.
Dover Publications, New York

\bibitem{burridge89}
Burridge, R., Chang, H. W.  (1989).
Multimode one-dimensional wave propagation in a highly discontinuous medium, Wave Motion 11, 231--249

\bibitem{white88}
Burridge, R., Papanicolaou, G., and White, B. (1988).
One-dimensional wave propagation
in a highly discontinuous medium, Wave Motion 10, 19--44

\bibitem{Buttiker1}
B\"{u}ttiker, M. (1988). Absence of backscattering in the quantum Hall effect in multiprobe conductors. Phys. Rev. B 38, 9375

\bibitem{Charlotte}
Charlotte, M, and Truskinovsky, L. (2012) Lattice dynamics from a continuum viewpoint. J Mech Phys Solids 60(8), 1508--1544

\bibitem{Chaudhuri}
Chaudhuri, A., Kundu, A., Roy, D., Dhar, A., Lebowitz, J. L., Spohn, H. (2010). Heat transport and phonon localization in mass-disordered harmonic crystals. Physical Review B, 81(6), 064301

\bibitem{clouet94}
Clouet, J.-F., Fouque, J.-P.  (1994).
Spreading of a pulse traveling in random media,
Ann. Appl. Probab. 4, 1083--1097

\bibitem{Dongre}
Dongre, B., Carrete, J., Katre, A., Mingo, N., Madsen, G. K. (2018). Resonant phonon scattering in semiconductors. Journal of Materials Chemistry C, 6(17), 4691-4697.

\bibitem{Deymier}
Deymier, P. and Dobrzynski, L. (2013). Discrete one-dimensional phononic and resonant crystals. In: Deymier P A (Ed.) {\em Acoustic metamaterials and phononic crystals}, vol. 173, Springer Series in Solid-State Sciences, Berlin Heidelberg: Springer, pp. 13--44
\bibitem{Elliott}
Elliott, R. J., Krumhansl, J. A., Leath, P. L. (1974). The theory and properties of randomly disordered crystals and related physical systems. Reviews of modern physics, 46(3), 465

\bibitem{Dhar}
Dhar, A., Dandekar, R. (2015). Heat transport and current fluctuations in harmonic crystals. Physica A: Statistical Mechanics and its Applications, 418, 49--64

\bibitem{Erdos}
Erd\"{o}s, P., Herndon, R. C. (1982). Theories of electrons in one-dimensional disordered systems. Advances in Physics, 31(2), 65--163

\bibitem{ethier}
Ethier, S. N., Kurtz, T. G. (1986). 
{\em Markov processes. Characterization and convergence},
Wiley, New York

\bibitem{book}
Fouque, J.-P., Garnier, J., Papanicolaou, G., Solna, K. (2007). {\em Wave propagation and time reversal in randomly layered media}. Springer, New York

\bibitem{Godin88}
Godin, T. J., Haydock, R. (1988). New method for calculation of quantum-mechanical transmittance applied to disordered wires. Physical Review B, 38(8), 5237

\bibitem{Sivak}
Gupta, D., Sivak, D. A. (2021). Heat fluctuations in a harmonic chain of active particles. Physical Review E, 104(2), 024605

\bibitem{Han}
Han, X. (2015). Asymptotic dynamics of stochastic lattice differential equations: a review. Continuous and Distributed Systems II: Theory and Applications, 121-136.

\bibitem{Itoh}
Itoh, T., Caloz, C. (2005). {\em Electromagnetic metamaterials: transmission line theory and microwave applications}. John Wiley \& Sons

\bibitem{Jara}
Jara, M., Komorowski, T., Olla, S. (2015). Superdiffusion of energy in a chain of harmonic oscillators with noise. Communications in Mathematical Physics, 339, 407--453

\bibitem{Katcho}
Katcho, N. A., Carrete, J., Li, W., Mingo, N. (2014). Effect of nitrogen and vacancy defects on the thermal conductivity of diamond: An ab initio Green's function approach. Physical Review B, 90(9), 094117.

\bibitem{Komorowski2}
Komorowski, T., Olla, S., Ryzhik, L. (2013). Asymptotics of the solutions of the stochastic lattice wave equation. Archive for Rational Mechanics and Analysis, 209, 455--494

\bibitem{Komorowski}
Komorowski, T., Stepien, L. (2012). Long time, large scale limit of the Wigner transform for a system of linear oscillators in one dimension. Journal of Statistical Physics, 148, 1--37

\bibitem{Tomasz}
Komorowski, T., Stepien, L. (2018). Kinetic limit for a harmonic chain with a conservative Ornstein-Uhlenbeck stochastic perturbation. Kinetic and Related Models, 11(2), 239--278

\bibitem{Kronig}
Kronig, R. D. L., Penney, W. G. (1931). Quantum mechanics of electrons in crystal lattices. Proceedings of the royal society of London. series A,
%containing papers of a mathematical and physical character,
130(814), 499--513

\bibitem{Kuo}
Kuo, F. (2006). Network analysis and synthesis. John Wiley \& Sons

\bibitem{Landauer1}
Landauer, R. (1957). Spatial variation of currents and fields due to localized scatterers in metallic conduction. IBM Journal of research and development, 1(3), 223--231

\bibitem{Lindsay}
Lindsay, L. (2016). Isotope scattering and phonon thermal conductivity in light atom compounds: LiH and LiF. Physical Review B, 94(17), 174304.

\bibitem{Muhlich}
M\"{u}hlich, U., Abali, B. E., dell'Isola, F. (2021) Commented translation of Erwin Schr\"{o}dinger’s paper
‘On the dynamics of elastically coupled point systems’ (Zur Dynamik elastisch gekoppelter Punktsysteme). 
Mathematics and Mechanics of Solids. 26(1):133--147

\bibitem{novikovbls}
Novikov, R., Sharma, B. L. (2023). 
Phase recovery from phaseless scattering data for discrete Schr\"{o}dinger operators. 
Inverse Problems 39(12), 125006.
%arXiv preprint arXiv:2307.06041.

\bibitem{oda}
O'Doherty, R. F., Anstey, N. A.  (1971).
Reflections on amplitudes, Geophysical Prospecting 19, 430--458

\bibitem{Ong}
Ong, Z. Y., Lee, C. H. (2016). Transport and localization in a topological phononic lattice with correlated disorder. Physical Review B, 94(13), 134203.

\bibitem{Payton}
Payton III, D. N., Visscher, W. M. (1967). Dynamics of disordered harmonic lattices. I. Normal-mode frequency spectra for randomly disordered isotopic binary lattices. Physical Review, 154(3), 802

\bibitem{Saxon}
Saxon, D. S., Hutner, R. A. (1949). Some electronic properties of a one-dimensional crystal model. Philips Research Reports (Netherlands) Changed to Philips J. Res., 4

\bibitem{Schrodinger}
Schr\"{o}dinger, E. (1914) Zur Dynamik elastisch gekoppelter Punktsysteme. Ann Phys 349(14): 916--934

\bibitem{Blsbif}
Sharma, B. L. (2016). Wave propagation in bifurcated waveguides of square lattice strips. SIAM Journal on Applied Mathematics, 76, 1355--1381

\bibitem{Blshexa} 
Sharma, B. L. (2018). Electronic transport across a junction between armchair graphene nanotube and zigzag nanoribbon: Transmission in an armchair nanotube without a zigzag half-line of dimers. The European Physical Journal B, 91, 1--25

\bibitem{Blshexa2}
Sharma, B. L. (2019). On electronic conductance of partially unzipped armchair nanotubes: further analysis. Eur. Phys. J. B 92, 1

\bibitem{Blsstep}
Sharma, B. L. (2020). Transmission of waves across atomic step discontinuities in discrete nanoribbon structures. Z. Angew. Math. Phys. 71, 73

\bibitem{Slater}
Slater, J. C. and Koster, G. F. (1954). 
{Simplified LCAO Method for the Periodic Potential Problem}, Phys. Rev. 94, 1498

\bibitem{Thebaud}
Thebaud, S., Berlijn, T., Lindsay, L. (2022). 
Perturbation theory and thermal transport in mass-disordered alloys: Insights from Green's function methods. Physical Review B, 105(13), 134202.


\bibitem{Xie}
Xie, G., Shen, Y., Wei, X., Yang, L., Xiao, H., Zhong, J., Zhang, G. (2014). 
A bond-order theory on the phonon scattering by vacancies in two-dimensional materials. Scientific reports, 4(1), 5085.

%\bibitem{Borcea}Borcea L., Garnier J., S{\o}lna K. (2021). Onset of energy equipartition among surface and body waves. Proc. R. Soc. A 477: 202007 75

%\bibitem{Hoop}de Hoop, M.V., Garnier, J. \& S{\o}lna, K. (2022). System of radiative transfer equations for coupled surface and body waves. Z. Angew. Math. Phys. 73, 177

\end{thebibliography}
\end{document}